\documentclass[11pt]{article}
\usepackage{amsmath,amssymb,epsf}
\usepackage[french]{babel}
\usepackage[applemac]{inputenc}
\advance \textwidth 3cm
\advance \textheight 3cm
\usepackage{amsfonts}
\usepackage{latexsym}
\newtheorem{theoreme}{Th\'eor\`eme}
\newtheorem{lemme}{Lemme}
\newtheorem{corollaire}{Corollaire}
\newtheorem{definition}{D\'efinition}

\newtheorem{prop}{Proprit}
\newtheorem{remarque}{Remarque}
\newenvironment{preuve}[1]{\par\noindent\underline{Preuve #1} :\quad}%
{\unskip\nobreak\hfil\penalty50\hskip2em\null\nobreak\hfil%
$\Box$\parfillskip0pt\par\medskip}

\newcommand{\trace}{\mathrm{Tr}}

\title{Op\'erateurs diff\'erentiels fractionnaires et matrices de Toeplitz.}

\author{ Philippe Rambour\thanks{Laboratoire de Mathmatiques dÕOrsay, Univ. Paris-Sud, CNRS, Universit Paris-Saclay, 91405 Orsay, 
tel : 01 69 15 57 28 ; fax 01 69 15 60 19
      \mbox{e-mail : philippe.rambour@math.u-psud.fr}}
       \and Abdellatif Seghier\thanks{Laboratoire de Mathmatiques dÕOrsay, Univ. Paris-Sud, CNRS, Universit Paris-Saclay, 91405 Orsay,  tel : 01 69 15 60 09 ; fax 01 69 15 72 34
       \mbox{ e-mail : abdellatif.seghier@math.u-psud.fr}}}
\begin{document}
\maketitle
  \renewcommand{\abstractname}{Rsum}
     \begin{abstract}
\textbf{ Op\'erateurs diff\'erentiels fractionnaire et matrices de Toeplitz.}\\
Dans ce travail on g\'en\'eralise \`a certains op\'erateurs fractionnaires la  m\'ethode utlis\'ee pour 
 pour inverser les op\'erateurs diff\'erentiels 
 $\frac{d^{2n}}{dx^{2n}}$ en inversant une matrice de Toeplitz. 
 L'int\'ert de ce travail est de montrer que l'on retrouve facilement par ce moyen les r\'esultats fournis par les m\'ethodes classiques d'analyse.
           \end{abstract}
          \renewcommand{\abstractname}{Abstract}
          \begin{abstract}
\textbf{ Fractional differential operators and Toeplitz matrices.}\\
In this work we generalize to few fractional differential operators the  method used to reverse differential operators
$\frac{d^{2n}}{dx^{2n}}$ by inverting a Toeplitz matrix. The interest of this work is to show that the method provides by the classical analytic methods of the analysis are easily founded by this means.
        \end{abstract}
   
%


\section{Introduction}
L'\'etude du calcul fractionnel est un sujet ancien puisqu'il a \'et\'e initi\'e il y a a peu pr\`es 300 ans par 
G.W Leibniz et L. Euler. Il est de nouveau d'actualit\'e depuis les ann\'ees 1970
et intervient dans la mod\`elisation de nombreux prob\`emes (voir \cite{TOBIAS1},\cite{TOBIAS2}).
Les d\'eriv\'ees fractionnaires qui vont nous occuper ici sont principalement les d\'eriv\'ees de Marchaud et de Grunwald-Letnikov . On sait que ces d\'eriv\'ees sont en fait reli\'ees entre elles et permettent \'egalement de retrouver les d\'eriv\'ees de Riemann-Liouville et
 de Caputo (voir \cite{Gorenflo1} pour les d\'eriv\'ees de Riemann-Liouville et de Caputo, 
 et \cite{Sam1} pour les d\'eriv\'ees fractionnaires en g\'en\'eral).

L'objectif de ce travail est d'introduire ces op\'erateurs, en particulier la d\'eriv\'ee fractionnaire inf\'erieure ou sup\'erieure de Marchaud sur un intervalle, au moyen 
d'une discr\'etisation de la fonction \`a  d\'eriver et d'une matrice de Toeplitz au symbole bien choisi. C'est \`a dire que ce travail s'inscrit dans la lign\'ee des travaux qui ont pour ambition d'associer un op\'erateur \`a une matrice de Toeplitz et 
l'inverse de cet op\'erateur \`a l'inverse de cette m\^{e}me matrice.
L'exemple consid\'er\'e ici est principalement celui qui consiste \`a r\'esoudre des \'equations diff\'erentielles 
fractionnaires sur un intervalle avec ou sans conditions de Dirichlet.  On r\'esoud des \'equations faisant intervenir des d\'eriv\'ees fractionnaires 
d'ordre $\alpha= 2p_{\alpha}+\alpha'$  avec $2 p_{\alpha} =[\alpha]$ et 
$0< \alpha<1$.  Si $p_{\alpha}=0$ on ne demande pas \`a la solution de v\'erifier de conditions aux bord de l'intervalle, et si $p_{\alpha}>0$ la solution doit v\'erifier 
$p_{\alpha}$ conditions initiales \`a chaque ext\'emit\'e de l'intervalle. Il est \`a noter que nous retrouvons les r\'esultats classiques \`a quelques variantes pr\`es. Ainsi les d\'eriv\'ees de Marchaud sup\'erieures ou inf\'erieures sont-elles d\'efinies pour des fonctions v\'erifiant des hypoth\`eses diff\'erentes de celles que l'on consid\`ere habituellement. Les \'enonc\'es \'etablis ici consid\'er\'ees ici portent sur des fonctions localement contractantes, o\`u dont les d\'eriv\'ees sont localement contractantes, sur un intervalle donn\'e, alors que dans le cas classique 
on consid\`ere plutt des fontions absolument continues ou 
de d\'eriv\'ees absolument continues, sur l'intervalle \'etudi\'e 
(voir \cite{Sam1}, et se rapporter aux th\'eor\`emes \ref{contractante} et \ref{diff}) (voir toujours \cite{Sam1} et 
les \'enonc\'es \ref{INTFRAC}, \ref{diffn} et 
\ref{GROSTHEOLOC} dans cet article). \\
Rappelons qu'une matrice de Toeplitz d'ordre $N$ de symbole $h$ est la matrice $T_{N} (h)$ d\'efinie par 
$\left( T_{N}(h) \right)_{(k+1,l+1} = \hat h (k-l)$ pour $0\le k,l\le N$ et o\`u $\hat h (s)$ d\'esigne le coefficient de Fourier d'ordre $s$.
Nous g\'en\'eralisons ici une m\'ethode bien connue dans le cas 
des \'equations diff\'erentielle du type 
$ \frac{d^{2n} f}{dt^{2n}}$ (\cite{SpSt},\cite{RRS}, \cite{RS04}, et aussi \cite{TOBIAS1} pour une autre approche) avec conditions initiales. 
Cette m\'ethode est bas\'ee sur l'inversion de certaines matrices de Toeplitz. 
Cette m\'ethode consiste \`a discr\'etiser l'\'equation au moyen de la matrice de Toeplitz 
d'ordre $N$ de symbole $f_{2n}$ d\'efini par 
$\theta \mapsto \vert 1- e^{i\theta}\vert ^{2n} = 
2^{2n} \vert \sin (\frac{\theta}{2})\vert ^{2n}$. 
Si $0<x<1$ et $\psi$ une fonction appartenant \`a $C^{2n}([0,1])$ on montre que la limite 
$\displaystyle{\lim_{N \rightarrow + \infty} \sum_{l=0}^{N}}
 \left(T_{N} (f_{2n})\right)_{k+1,l+1} \psi (\frac{l}{N}),$ avec $k=[Nx]$, n'est rien d'autre, \`a un coefficient pr\`es, que $ (-1)^{n} \frac{d^{2n} f\psi}{dt^{2n}}$. L'inversion de la matrice de Toeplitz $T_{N}(f_{n})$ fournit alors un noyau de Green associ\'e aux conditions initiales 
 $f^{(0)}(0)=f^{(1)}(0)=\cdots=f^{(n-1)}(0)=0$ et $f^{(0)}(1)=f^{(1)}(1)=\cdots=f^{(n-1)}(1)=0$.\\
 Ici nous consid\'erons pour $\alpha>0$ la matrice de Toeplitz d'ordre $N$ de symbole 
 $\varphi_{\alpha}= \lim_{R\rightarrow 1^{-}} \varphi_{\alpha,R}$ o\`u $\varphi_{\alpha,R}$ 
  est la fonction d\'efinie sur $]-\pi,\pi[$ par 
 $ \theta\mapsto (1-Re^{i\theta})^{\alpha}(1+R e^{-i\theta})^{\alpha}.$
 Pour $f$ une fonction d\'efinie sur $]0,1[$ et $0< x< 1$ on \'etudie alors la limite
 $$ \lim_{N\rightarrow +\infty} N ^{\alpha} 
 \left(\sum_{l=0}^{N}T_{N} \left(\varphi_{\alpha}\right)_{k+1,l+1} \left(X_{N}\right)_{l}
 \right) =\left( D_{\alpha} (f)\right) (x), \quad \mathrm{avec} \quad 
k=[Nx],$$ 
en posant $(X_{N})_{l}= f(\frac{l}{N})$ si $l\neq 0,1$,
$(X_{N})_{0}=(X_{N})_{1}=0$.
Nous posons alors les d\'efinitions suivantes.
\begin{definition}
On dira qu'une fonction $f$ est localement contractante sur un intervalle $]a,b[$ si pour tout intervalle $[\delta_{1}, \delta_{2}]$ contenu dans 
$[a,b]$ il existe un r\'eel $K$ 
v\'erifiant, pour tout $z,z'$ dans $[\delta_{1}, \delta_{2}]$ :
$$ \vert f(z) -f(z')\vert \le K \vert z-z'\vert. $$
\end{definition}
\begin{definition}
Pour tout intervalle $]a,b[$ et tout r\'eel $\gamma\in [0,1[$ on notera 
$L^{1}_{1-\alpha} (]a,b[)$ l'ensemble des fonctions $f$ d\'efinies sur $]a,b[$
v\'erifiant pour un r\'eel $\gamma\in [0, 1-\alpha]$ $f(x) =O(x^{-\gamma}), f(1-x) = O(x^{-\gamma})$ quand $x\to 0$ par valeur positives .
\end{definition}
On consid\`ere alors les fonctions localement contractantes sur $]0,1[$ et dans 
$L^{1}_{1-\alpha} (]a,b[)$ pour $\gamma\in [0, 1-\alpha]$.

On montre (voir le th\'eor\`eme \ref{contractante}) que si $\alpha\in ]0,1[$ cette limite existe pour de telles fonctions et vaut est, pour $x\in]0,1[ $ ou pour $x\in [0,1]$ selon les hypoth\`eses,

$$
\frac{2^\alpha}{\Gamma(-\alpha)}\left(
 \int_{0}^{x} \vert x-t\vert^{-\alpha-1} \left(f(t)- f(x)\right)dt - f(x) \left(\frac{(x)^{-\alpha}}
 {\alpha}\right)\right)
 $$
c'est \`a dire que nous retrouvons la d\'eriv\'ee fractionnaire inf\'erieure de Marchaud d'ordre $\alpha$ sur $(0,1)$.
On peut d'ailleurs v\'erifier que si l'on choisit comme symbole la fonction 
$\tilde \varphi_{\alpha}= \lim_{R\rightarrow 1^{-}} \tilde \varphi_{\alpha,R}$ o\`u 
$\tilde \varphi_{\alpha,R}$ 
  est la fonction d\'efinie sur $]-\pi,\pi[$ par 
 $ \theta\mapsto (1+Re^{i\theta})^{\alpha}(1-R e^{-i\theta})^{\alpha}$
cette m\^{e}me limite nous donne  la d\'eriv\'ee fractionnaire sup\'erieure de Marchaud d'ordre $\alpha$ sur $(0,1)$. 
Consid\`erons la d\'eriv\'ee fractionnaire de 
Gr\"{u}nwald-Letnikov calcul\'ee en un r\'eel $x$ d'un intervalle $[a,b]$ et donn\'ee pour une fonction $f$ d\'efinie sur $[a,b]$ par 
$$ \left((D_{G L}^{\alpha})f \right)(x) = \lim_{\Delta x \rightarrow 0} 
\frac{1}{(\Delta x)^{\alpha}}\sum_{n=0}^{\frac{x-a}{\Delta x}}
(-1)^{n} \frac{\Gamma(\alpha+1)}{\Gamma(n+1)\Gamma(\alpha-n+1)} f(x-nx).$$
Si l'on choisit $\Delta x = \frac{1}{N} $,
$[a,b]=[0,1]$ alors 
$$ \left( D_{G L}^{\alpha} \right)(x) = \lim_{N\rightarrow +\infty}  
N^{\alpha}\sum_{l=0}^{k}
(-1)^{l} \frac{\Gamma(\alpha+1)}{\Gamma(k-l+1)\Gamma(\alpha-(k-l)+1)} f(\frac{l}{N})$$
avec $k=[Nx]$.
Un calcul identique \`a celui effectu\'e dans la d\'emonstration du th\'eor\`eme 
\ref{contractante} permet d'obtenir, en choisissant judicieusement le symbole de la matrice deToeplitz utilis\'ee, que pour une fonction contractante sur 
$[0,1]$ on a  $ \left(D_{G L}^{\alpha} \right)(x) = \left( D_{\alpha} (f)\right) (x)$.\\
En adaptant \`a la dimension un des formules d'inversion des matrices de Toeplitz \`a blocs \'etablies dans \cite{RRQuebl} on peut obtenir les coefficients de la matrice $T_{N}^{-1} (\varphi_{\alpha,R})$. On peut alors d\'efinir une 
int\'egration fractionnaire not\'ee $J_{\alpha}$ pour une fonction 
$f$ d\'efinie sur $[0,1]$. Pour $x\in [0,1]$ $\left(J_{\alpha} (f)\right) (x)$ sera d\'efinie si la limite 
quand $R$ tend vers $1$ par valeurs inf\'erieurs et quand $N$ tend vers plus l'infini
de $\displaystyle{ \sum_{l=0}^{N} \left(T_{N}^{-1} (\varphi_{\alpha,R})\right)_{k+1,l+1}}$
existe et est finie.
En particulier pour $\alpha\in ]0,1[$ et $f$ une fonction localement contractante sur $]0,1[$ et dans 
$L^{1}_{1-\alpha}([a,b])$ on obtient que pour $x\in ]0,1[$ $\left(J_{\alpha} (f)\right) (x)$ existe et vaut  (propri\'et\'e 
\ref{INTFRAC})
$$\left(J_{\alpha} (f)\right) (x) = \frac{1}{2^{\alpha}\Gamma (\alpha)}\int_{0}^{x} \frac{f(t)}{(x-t)^{1-\alpha}}dt.$$
L'existence de $J_{\alpha}$ nous permet  de r\'esoudre des \'equations diff\'erentielles fractionnaires. Par exemple pour $\alpha\in ]0,1[$ on obtient 
que pour une fonction $\psi$ contractante sur $]0,1[$ la fonction $h= J_{\alpha}(\psi)$
est d\'efinie sur $]0,1[$ et est solution  pour tout $x\in ]0,1[$ de l'\'equation
$D_{\alpha} (y) = \psi$ (voir la propri\'et\'e \ref{invgauche} et le corollaire \ref{equdiff1}. \\
On \'etend ensuite ce r\'esultat \`a $D_{\alpha}$ pour $\alpha\ge1$. Soit $\alpha$ un r\'eelt etn un entier positif. SI $f$ est une fonction n fois d\'erivable   dans $]0,1[$ 
 v\'erifiant $f(0)=f^{(1)}(0)=\cdots= f^{(n-1)}(0)=0$,
 $ f(1)= f^{(1)}(1)= \cdots=f^{(n-1)}(1)=0$ et $f^{(n)}$ contractante.
 Alors  $f$ est $D^{\alpha+n}$ d\'erivable si et seulement si $f^{(n)}$ est
 $\alpha$ d\'erivable et de plus pour tout r\'eel $x$ dans $]0,1[$ 
 $$ D^{\alpha+n} (f) (x) = (2)^nD^{\alpha} (f^{(n)}) (x).$$
 
Si $n=2p$ $J_{n}(f) $ est d\'efini \`a partir du noyau de Green 
d\'efini  pour $0\le x,y\le 1$ par la formule 
$$G_{p}(x,y) = \frac{x^py^p} {\left((p-1)!\right)^2} 
  \int_{\max (x,y)}^1\frac{(t-x)^{p-1} (t-y)^{p-1}}{t^{2p}} dt
  \quad \mathrm{si} \quad (x,y)\neq (0,0)$$
  et $G(0,0)=0$.
 
  On rappelle qu'il a \'et\'e \'etabli dans (\cite{RS04}) que 
  $G_{p}$ est le noyau de Green associ\'e \`a l'op\'erateur 
  diff\'erentielle $ (-1)^p d^{2p}/dx^{2p}$ sur $[0,1]$ associ\'e 
  aux conditions initiales 
  $$ f^{(0)} (0) = f^{(1)} (0)= \cdots = f^{(p-1)} (0)=0,
  \quad 
   f^{(0)} (1) = f^{(1)} (1)= \cdots = f^{(p-1)} (1)=0.$$ 
   On pose alors 
 $ J_{n} (f) = (-1)^p \int_{0}^1G_{p}(x,t) f(t) dt$, puis 
 o\`u $\tilde J_{\alpha} = J_{n_{\alpha}}\circ J_{\alpha-n_{\alpha}}$.\\
 On peut \`a ce moment l\`a  utiliser les noyau de Green obtenus dans \cite{RS04}
 r\'esoudre des \'equations diff\'erentielles du type $D_{\alpha} (y) = \psi$ avec 
 conditions de Dirichlet comme dans le th\'eor\`eme :
  \begin{theoreme}\label{GROSTHEO}
  On consid\`ere un r\'eel positif $\alpha$, de partie enti\`ere $n_{\alpha}$.
 on suppose que $n_{\alpha} =2p_{\alpha}$ est un entier pair  et
  on se donne une fonction  $\psi$ localement contractante sur $]0,1[$ et dans $L^1_{1-\alpha}([0,1])$. Alors 
  la fonction $\tilde J_{\alpha} (\psi)$ est solution sur $]0,1[$
  de l'\'equation diff\'erentielle 
  $D_{\alpha}(y)=\psi$
  avec les conditions initiales 
  $ y(0)=y^{(1)}(0)= \cdots=y^{(p_{\alpha}-1)}(0)=0$\\
  et 
  $ y(1)=y^{(1)}(1)= \cdots=y^{(p_{\alpha}-1)}(1)=0$.
 \end{theoreme}

 On peut ensuite par homoth\'etie d\'efinir les d\'eriv\'ees $D_{\alpha,a,b}$ 
 et les primitives $\tilde J_{\alpha,a,b} $ sur tout intervalle 
 $[a,b]$ et r\'esoudre des \'equations diff\'erentielles \`a conditions initiales 
 sur un intervalle $[a,b]$ gr\^{a}ce \`a un th\'eor\`eme analogue 
 au t\'h\'eor\`eme \ref{GROSTHEO}. 
 Pour $0<\alpha<1$ par un passage \`a la limite si,   $\displaystyle{ \lim_{a\rightarrow -\infty, b \rightarrow + \infty} 
\frac{\left(D_{\alpha,a,b} (f)\right)(x) }{(b-a)^{-\alpha}}}$ existe et est fini.
on d\'efinit la d\'eriv\'ee 
$$\left( D_{\alpha,\infty} (f)\right) =  \frac{2^\alpha }{\Gamma (-\alpha)} 
\int_{-\infty}^{x} (x-u)^{-\alpha-1} \left ( f(u) - f (x)\right) du.$$
On  d\'efinit de m\^{e}me $ \left( J_{\alpha,\infty} (f)\right) (x) = \displaystyle{ \lim_{a\rightarrow -\infty, b \rightarrow + \infty} 
\frac{\left(J_{\alpha,a,b} (f)\right)(x) }{(b-a)^{\alpha}}}$ 
 si cette limite existe et est finie.
 On remarque que $D_{\alpha,\infty}$ est une d\'eriv\'ee de Marchaud.
Dans ce cas si $\psi$ est une fonction localement contractante, \`a support compact, (donc en fait contractante et int\'egrable) on peut \'egalement r\'esoudre sur 
 $\mathbb R$ l'\'equation 
 diff\'erentielle $\left(D_{\alpha,\infty} (y)\right (x)=\psi(x)$ en posant 
 $y=J_{\alpha,\infty}$.\\
 En conclusion ce travail met en \'evidence le lien profond qui existe entre certaines matrices de Toeplitz et des op\'erateurs diff\'erentiels d'une part, et entre l'inverse de ces matrices et l'inverse de ces op\'erateurs d'autre part.  Les formules d'inversion des matrices de Toeplitz peuvent alors permettre d'inverser ces op\'erateurs.  
 \section{La d\'erivation fractionnaire $D_{\alpha}$}
\subsection{D\'efinition de la fonction $\varphi_{\alpha}$}
Rappelons que nous notons par $\chi$ la fonction $\theta\mapsto e^{i\theta}$.
Pour $\alpha>0$, $R\in ]0,1[$ et $\theta\in ]-\pi,\pi]$ nous d\'efinissons les fonctions $\varphi_{\alpha,R}$ et $\varphi_{\alpha}$ par 
$\varphi_{\alpha,R}  = 
(1-R\chi)^{\alpha} (1+R \bar \chi)^{\alpha}= (1-R^{2}-R\chi+R\bar \chi)^{\alpha}$.
et $\varphi_{\alpha}= \displaystyle{\lim_{R\rightarrow1_{-}} \varphi_{\alpha,R}}$\\
Nous allons maintenant d\'eterminer l'asymptotique des coefficients de Fourier 
de la fonction 
$ \varphi_{\alpha}$
\begin{theoreme} \label{contraction}
Si $\alpha$ d\'esigne un r\'eel strictement sup\'erieur \`a 
$-\frac{1}{2}$ non entier, on a
 pour tout entier $n$ suffisamment grand, 
 \[
\left\{
\begin{array}{cccc}
   \widehat {\varphi_{\alpha}}(n)  &  = &   n^{-\alpha-1} 
   \frac{2^{\alpha}}{\Gamma (-\alpha)} + o(n^{-\alpha-1})   \quad &\mathrm{si} \quad n > 0\\
  \widehat {\varphi_{\alpha}}(n)  &  = & (-1)^{n} (- n)^{-\alpha-1} 
   \frac{2^{\alpha}}{\Gamma (-\alpha)}+ +o(n^{-\alpha-1})  \quad &\mathrm{si} \quad n < 0. \end{array}
\right.
\]
\end{theoreme}
\begin{preuve}{}
Rappelons tout d'abord que pour tout r\'eel non entier $\alpha$
et tout $R\in [0,1[ $ nous avons 
$$ (1+R\chi)^{-\alpha} = \sum_{u\ge 0}
\beta_{u}^{(\alpha)} R^u \chi^u$$
avec (voir \cite{ZYG2})
\begin{equation} \label{AS}
\beta_{u}^{(\alpha)} = \frac{u^{\alpha-1}}{\Gamma (\alpha)}
+ O (u^{\alpha-2}).
\end{equation}
Posons $\delta _{n}^{(\alpha,R)}$ le coefficient de Fourier d'ordre $n$ de 
$\varphi_{\alpha,R},$ et  $\delta _{n}^{(\alpha)}$ celui de 
$\varphi_{\alpha}.$\\
Nous avons 
\begin{equation}\label{EGALITE1}
\delta _{n}^{(\alpha,R)} = \sum_{u\ge 0} \beta_{u}^{(-\alpha)} 
\beta_{n+u}^{(-\alpha)} (-1)^{u} R^{n+2u} 
\end{equation}
si $n$ est un entier positif et
\begin{equation}\label{EGALITE2}
\delta _{n}^{(\alpha,R)} = \sum_{u\ge 0} \beta_{u}^{(-\alpha)} 
\beta_{-n+u}^{(-\alpha)} (-1)^{-n+u} R^{-n+2u} 
\end{equation}
si $n$ est un entier n\'egatif.\\
Pour $n$ un entier positif suffisamment grand et 
$\gamma$ un r\'eel dans $]0,1[$ 
\'ecrivons
$$ \delta _{n}^{(\alpha,R)} = \sum_{u=0} ^{n^\gamma}\beta_{u}^{(-\alpha)} 
\beta_{n+u}^{(-\alpha)} (-1)^{u} R^{n+2u}
+ \sum_{u >n^\gamma} \beta_{u}^{(-\alpha)} 
\beta_{n+u}^{(-\alpha)} (-1)^{u} R^{n+2u}.$$
En \'ecrivant pour $u$ et $n$ assez grand le terme  
$ \beta_{u}^{(-\alpha)} 
\beta_{n+u}^{(-\alpha)} (-1)^{u} R^{n+2u}$
sous la forme 
$ a_{u} b_{u} + O(u^{-\alpha-2} n^{-\alpha-1})$ avec 
$a_{u}= u^{-\alpha-1} (n+u)^{-\alpha-1}$,
$ b_{u}= (-1)^{u} R^{n+2u}$
il est facile de v\'erifier, en utilisant une sommation d'Abel, que 
\begin{equation} \label{RESTA1}
\Bigl \vert 
\sum_{u >n^\gamma} \beta_{u}^{(-\alpha)} 
\beta_{n+u}^{(-\alpha)} (-1)^{u} R^{n+2u}\Bigr \vert 
\le K_{1} \left(n^{\gamma}\right)^{-\alpha-1} n^{-\alpha-1}.
 \end{equation}
Nous avons d'autre part 
$$\sum_{u=0} ^{n^\gamma}\beta_{u}^{(-\alpha)} 
\beta_{n+u}^{(-\alpha)} (-1)^{u} R^{n+2u} =
\frac{ n^{-\alpha-1} R^n}{\Gamma (-\alpha)} 
\left( \sum_{u=0} ^{n^\gamma}\beta_{u}^{(-\alpha)} 
(-1)^u R^{2u} + \tau( \gamma)\right) $$
avec 
\begin{equation} \label{RESTA2}
 \vert  \tau( \gamma) \vert \le K_{3} \left( \displaystyle
{\sum_{u=0} ^{n^\gamma}\frac{u}{n} }\right) 
\le K_{4}\left( n^{2\gamma-1}\right) 
\end{equation}
Enfin nous obtenons, en utilisant les propri\'et\'es des s\'eries altern\'ees,  
\begin{equation} \label{RESTA3}
\sum_{u=0} ^{n^\gamma}\beta_{u}^{(-\alpha)} 
(-1)^u R^{2u} = (1+R^2)^\alpha + K_{5} n^{(-\alpha -1)\gamma}
\end{equation}
On peut d'autre part v\'erifier que les constantes $K_{i}, 1\le i \le 5$ qui apparaissent dans 
les \'equations (\ref{RESTA1})(\ref{RESTA2})(\ref{RESTA3})
sont ind\'ependantes de $R$ et de $n$. 
Nous pouvons donc affirmer que pour de entier $n$ 
suffisamment grand et $\gamma$ bien choisi 
$\delta_{n} ^{(\alpha,R)} = n^{-\alpha-1} 
\frac{ (1+R^2)^\alpha R^n}{\Gamma (-\alpha)} 
+o(n^{-\alpha-1})$
uniform\'ement par rapport \`a $R$.
Soit en faisant tendre $R$ vers 1 :
$$\delta_{n} ^{(\alpha)} =  n^{-\alpha-1} \frac{ 2^\alpha }{\Gamma (-\alpha)} +o(n^{-\alpha-1}).$$
De la m\^{e}me mani\`ere on obtient le th\'eor\`eme
pour $n < 0$ \`a partir de (\ref {EGALITE2}).
\end{preuve}
\subsection{D\'efinition de la d\'erivation $D_{\alpha}$}
\begin{definition}
Soient $\alpha$ est un r\'eel strictement sup\'erieur \`a z\'ero 
et  $f$ une fonction d\'efinie sur $[0,1]$. On dira que $f$ est $\alpha$ d\'erivable en $x\in [0,1]$ 
 si et seulement si il existe un r\'eel, not\'e $\left( D_{\alpha} (f)\right) (x)$ tel que 
 $$ \lim_{N\rightarrow +\infty} N ^{\alpha} 
 \left(\sum_{l=0}^{N}T_{N} \left(\varphi_{\alpha}\right)_{k+1,l+1} \left(X_{N}\right)_{l}
 \right) =\left( D_{\alpha} (f)\right) (x), \quad \mathrm{avec} \quad 
k=[Nx].
 $$
\end{definition}
\begin{theoreme}\label{contractante}
Soit $\alpha\in ]0,1[$ et $f$ une fonction localement contractante d\'efinie sur $]0,1[$, appartenant \`a 
$L_{1-\alpha}^1(]a,b[)$.  
 Alors $\left( D_{\alpha} (f)\right) (x)$ existe et vaut 
$$ 
 \left( D_{\alpha} (f)\right) (x) = \frac{2^\alpha}{\Gamma(-\alpha)}\left(
 \int_{0}^{x} \vert x-t\vert^{-\alpha-1} \left(f(t)- f(x)\right)dt - f(x) \left(\frac{(x)^{-\alpha}}
 {\alpha}\right)\right)
 $$
 pour tout $x\in ]0,1[$.
 \end{theoreme}
 \begin{remarque}
 Si les quantit\'es $f(0)$ ou $f(1)$ ne sont pas d\'efinies, on pose $f(0)$ ou $f(1)=0$.
 \end{remarque} 
 \begin{preuve}{}
 Dans la suite de cette d\'emonstration nous poserons 
 $C_{\alpha}= \frac{2^\alpha}{\Gamma(-\alpha)}$.
 Pour $N$ fix\'e posons $k$ tel que $[Nx]=k$ et $\epsilon$ un r\'eel positif suffisamment petit.
 On a 
 \begin{align*}   
\sum_{l=0}^{N}T_{N} \left(\varphi_{\alpha}\right)_{k+1,l+1} \left(X_{N}\right)_{l} 
&=
 \sum_{l=0}^{k-[N\epsilon]-1}T_{N} \left(\varphi_{\alpha}\right)_{k+1,l+1} \left(X_{N}\right)_{l} +\sum_{l=k-[N\epsilon]}^{k+[N\epsilon]}T_{N} \left(\varphi_{\alpha}\right)_{k+1,l+1} \left(X_{N}\right)_{l}
 \\
 &+ \sum_{l=k+[N\epsilon]+1}^{N}T_{N} \left(\varphi_{\alpha}\right)_{k+1,l+1} \left(X_{N}\right)_{l}.
 \end{align*}
  Ce qui peut s'\'ecrire aussi 
    \begin{align*} 
 \sum_{l=0}^{N}T_{N}\left(\varphi_{\alpha}\right)_{k+1,l+1} 
 \left(X_{N}\right)_{l} &=
  \sum_{l=0}^{k-[N\epsilon]-1} \widehat {\varphi_{\alpha}} 
  (k-l)
  f(\frac{l}{N}) +\sum_{l=k+[N\epsilon]+1}^{N}\widehat {\varphi_{\alpha}} (k-l)  
  f(\frac{l}{N}) \\
 & - \left( \sum_{l < k-[N\epsilon]} \widehat {\varphi_{\alpha}} 
 (k-l)
    +\sum_{l>k+[N\epsilon]} \widehat {\varphi_{\alpha}} 
    (k-l)\right)
    f(\frac{k}{N}) +R_{\epsilon,N}
   \end{align*}
avec, en utilisant la convergence de la s\'erie 
$\displaystyle{\sum_{u\in \mathbb Z}  \widehat {\varphi_{\alpha}}(u)}$, et la locale convergence,
$$ R_{\epsilon,N} = \sum_{l=k+[N\epsilon]}^{k-[N\epsilon]}
\widehat{\varphi_{\alpha}}(k-l) 
\left( f(\frac{l}{N})-f(\frac{k}{N})\right).
$$
\begin{equation} \label{restemajore} 
\vert R_{\epsilon,N}\vert \le \epsilon \sum_{l = k-[N\epsilon]+1}^{l=k+[N\epsilon]-1}
\vert \widehat {\varphi_{\alpha}} (k-l)\vert 
    \le \epsilon  M 
    \end{equation} 
  ($M$ est un r\'eel ind\'ependant de $N$).
  Ce qui devient 
  \begin{align*} 
\sum_{l=0}^{N}T_{N}\left(\varphi_{\alpha}\right)_{k+1,l+1} \left(X_{N}\right)_{l} 
&=\sum_{l=0}^{k-[N\epsilon]-1} \widehat {\varphi_{\alpha}} 
(k-l)
 \left( f(\frac{l}{N}) -f(\frac{k}{N})\right) 
 \\& +\sum_{l=k+[N\epsilon]+1}^{N}\widehat {\varphi_{\alpha}} 
 (k-l)  
\left( f(\frac{l}{N} -f(\frac{k}{N})\right) \\ &
 - \left( \sum_{l < 0} \widehat {\varphi_{\alpha}} 
 (k-l)
  +\sum_{l>N} \widehat {\varphi_{\alpha}} (k-l)\right)
  f(\frac{k}{N}) +R_{\epsilon,N}
  \end{align*}
  On obtient, en utilisant le th\'eor\`eme \ref{contraction}  
 \begin{align*}    \sum_{l=k+[N\epsilon]+1}^{N } \widehat {\varphi_{\alpha}} (k-l)
 & \left( f(\frac{l}{N}) -f(\frac{k}{N})\right) 
 = \\
 &C_{\alpha}  \sum_{l=k+[N\epsilon]+1}^{N} (-1)^{k-l} 
 (l-k)^{-\alpha-1}   \left( f(\frac{l}{N}) -f(\frac{k}{N})\right) +o(N^{-\alpha}). 
\end{align*} 
Une sommation d'Abel donne que 
\begin{equation}\label{MAJO1}
\Bigl \vert \sum_{l=k+[N\epsilon]+1}^{N} (-1)^{k-l}
 (l-k)^{-\alpha-1}\left( f(\frac{l}{N}) -f(\frac{k}{N})\right)\Bigr \vert 
 =O\left( (N) ^{-\alpha-1} \epsilon ^{-\alpha}\right).
 \end{equation}
 En effet il vient, avec la formule d'Euler et Mac-Laurin,
 $$
 \Bigl \vert \sum_{l=k+[N\epsilon]+1}^{N} (-1)^{k-l}
 (l-k)^{-\alpha-1}\left( f(\frac{l}{N}) -f(\frac{k}{N})\right)\Bigr \vert
 \le \vert C_{\alpha}\vert \left( S_{1} + S_{2} +S_{3} \right)
 $$
 avec, en utilisant les hypoth\`eses faites sur $f$  
 $$S_{1}= (N\epsilon)^{-\alpha-1} \epsilon + k^{-\alpha-1}
 \vert f(1) - f(\frac{k}{N}) \vert = O\left( N^{-\alpha-1} 
 \epsilon ^{-\alpha}\right),$$
\begin{align*}
 S_{2 } &= \sum_{l=k+[N\epsilon]+1}^{N} 
 \vert (l-k)^{-\alpha-1} - (l-k-1)^{-\alpha-1}\vert 
 \vert  f(\frac{l}{N}) -f(\frac{k}{N})\vert \\ 
 &\le O\left( \sum_{l=k+[N\epsilon]+1}^{N} (l-k-1)^{-\alpha-2} \right)\\
 &= O\left( (N\epsilon)^{-\alpha-1} \right),
 \end{align*}
 et enfin 
 $$ 
 S_{3}= \sum_{l=k+[N\epsilon]+1}^{N} 
 \vert (l-k)^{-\alpha-1}\vert  
 \vert  f(\frac{l}{N}) -f(\frac{l-1}{N})\vert.$$
 On utilise la d\'ecomposition  
 \begin{align*}
  S_{3 } &=\sum_{l=k+[N\epsilon]+1}^{N- [N\epsilon]} 
 \vert (l-k)^{-\alpha-1}\vert  
 \vert  f(\frac{l}{N}) -f(\frac{l-1}{N})\vert\\
 &+ \sum_{l=N-[N\epsilon]+1}^{N} 
 \vert (l-k)^{-\alpha-1}\vert  
 \vert  f(\frac{l}{N}) -f(\frac{l-1}{N})\vert
 \end{align*}
 Soit $\epsilon_{0}$ un r\'eel positif tendant vers z\'ero. Si $K_{\epsilon_{0}}$ est tel que pour tous les r\'eel $x$ et $x'$ dans 
 $[x+\epsilon, 1-\epsilon_{0}]$ on a $\vert f(x) - f(x') \vert \le K_{\epsilon_{0}} \vert x - x'\vert $
 alors 
\begin{align*}
  \sum_{l=k+[N\epsilon]+1}^{N- [N\epsilon_{0}]} 
 \vert (l-k)^{-\alpha-1}\vert  
 \vert  f(\frac{l}{N}) -f(\frac{l-1}{N})\vert\le 
  &\le K_{\epsilon_{0}}  \sum_{l=k+[N\epsilon]+1}^{N-[N\epsilon]} (k-l)^{-\alpha-1} \frac{1}{N}  \\
 &= O\left( N^{-\alpha-1} \epsilon^{-\alpha} \right).
 \end{align*}
 D'autre part l'hypoth\`ese 
 $f\in L^{1}_{1-\alpha}([0,1])$,  implique l'existence d'un r\'eel $\gamma \in [0, 1-\alpha]$ tel que 
 $$ \sum_{l=N-[N\epsilon_{0}]+1}^{N} 
 \vert (l-k)^{-\alpha-1}\vert  
 \vert  f(\frac{l}{N}) -f(\frac{l-1}{N})\vert = O\left(N^{\gamma} (N-k)^{\alpha-1}\epsilon_{0} \right) = o(1)$$
 On remarque enfin que si $\alpha\in ]0, 1[ $ et 
$\epsilon= N^{-\tau}$ on peut choisir $\tau$ 
de mani\`ere \`a ce que $ \frac{1}{\alpha+1}>\tau>\alpha$ 
ce qui implique 
 \begin{equation}\label{formule1}
 \sum_{l=k+[N\epsilon]+1}^{N} \widehat {\varphi_{\alpha}} (k-l)  \left( f(\frac{l}{N}) - f(\frac{k}{N})\right)
= o(N^{-\alpha_{0}})
\end{equation}
pour tout r\'eel $\alpha_{0} \in ]0, \alpha[$.
et 
$R_{\epsilon,N} = o(N^{-\alpha}).$
D'autre part 
\begin{equation} \label{formule2}
  \sum_{l=0}^{k-[N\epsilon]-1}\widehat {\varphi_{\alpha}} 
  (k-l)  
 \left( f(\frac{l}{N}) -f(\frac{k}{N})\right)
  = C_{\alpha} N^{-\alpha}  \int_{0} ^{x-\epsilon}\vert x-t\vert ^{-\alpha-1} \left(f(t) -f(x) \right) dt +o(N^{-\alpha}). 
 \end{equation}
 et aussi 
 \begin{equation}\label{formule3}
 \left( \sum_{l < 0} \widehat {\varphi_{\alpha}} (k-l)
    +\sum_{l>N} \widehat {\varphi_{\alpha}} (k-l)\right)
    f(\frac{k}{N})  = C_{\alpha} N^{-\alpha}  \frac{x^{-\alpha}}{\alpha}f(x)
   +o(N^{-\alpha}).
    \end{equation}
   Nous avons donc obtenu
 $$ N^{\alpha} \sum_{l=1}^{N}T_{N}
 \left(\varphi_{\alpha}\right)_{k+1,l+1} \left(X_{N}\right)_{l} 
 = 
C_{\alpha} \left( \int_{0} ^{x-\epsilon} \vert x-t\vert ^{-\alpha-1} \left(f(t) - f(x) \right) dt-  \frac{x^{-\alpha}}{\alpha} f(x) \right)
+o(1).
$$
 \end{preuve}
Nous avons de plus 
 \begin{prop} \label{propzero}
 Soit $0<\alpha<1$ et $f\in L^{1}_{1-\alpha}([0,1[)$ une fonction localement contractante sur $[0,1[$,
 v\'erifiant $f(0)=0$.
 Alors $\left( D_{\alpha}(f)\right)(0)$ existe et vaut
 z\'ero
  \end{prop}
 \begin{preuve}{}
Soit  $\epsilon_{1}\in ]0,1[$.
On a 
 $$ \sum_{l=0}^{N}\widehat{\varphi_{\alpha} } (-l) f(\frac{l}{N}) 
 = \sum_{l=0}^{ [N\epsilon_{1}]}\widehat{\varphi_{\alpha}} (-l)  
 \left(f(\frac{l}{N})-f(0)\right)
+\sum_{l= [N\epsilon_{1}]+1}^{N}\widehat{\varphi_{\alpha} }(-l) 
f(\frac{l}{N}) .$$
Gr\^{a}ce \`a l'hypoth\`ese contractante 
$$  \sum_{l=0}^{ [N\epsilon_{1}]}\widehat{\varphi_{\alpha}} (-l) 
 \left(f(\frac{l}{N})-f(0)\right)
\le  O( \epsilon_{1})=o(N^{-\alpha}) 
,$$
si on choisit $\epsilon_{1}= N^{-\gamma_{1}}$ avec $1>\gamma_{1}>\alpha$.
D'autre part en utilisant la d\'emonstration du th\'eor\`eme \ref{contraction} il vient 
$$  \sum_{l= [N\epsilon_{1}]+1}^{N}\widehat{\varphi_{\alpha} }(-l) 
f(\frac{l}{N})  = O(N^{-\alpha-1}\epsilon_{1}^{-\alpha})$$
ce qui implique qu'il faut finalement choisir $\alpha<\gamma_{1}\leq\frac{1}{\alpha}$.
La propri\'et\'e est alors d\'emontr\'ee.
\end{preuve}
 \begin{prop} \label{propzero2}
  Soit $<\alpha<1$ et  $f\in ]0,1]$ est une fonction localement  contractante sur $]0,1]$, 
  et appartenant \`a $L^{1}_{\gamma}(]0,1])$, avec  $0<\gamma<1-\alpha$, v\'erifiant 
  $f(1)=0$.
  Alors 
 $\left( D_{\alpha}(f)\right)(1)$ existe et vaut 
 $ \frac{2^{\alpha}}{\Gamma(-\alpha)} 
 \left( \int_{0}^{1} t^{-\alpha-1} f(t) dt \right).$
  \end{prop}
\begin{preuve}{} 
En posant $[N\epsilon_{2}]= N^{-\gamma} $ avec $1>\gamma>0$ 
on \'ecrit 
$$
 \sum_{l=0}^{N}\widehat{\varphi_{\alpha} }(N-l) f(\frac{l}{N}) 
 = \sum_{l=0}^{N- [N\epsilon_{1}]}
 \widehat{\varphi_{\alpha}} (N-l) f(\frac{l}{N})
+\sum_{l= N- [N\epsilon_{1}]-1}^{N}\widehat{\varphi_{\alpha}}  (N-l)  \left(f(\frac{l}{N})-f(1)\right).
$$
On raisonne ensuite comme pour la d\'emonstration du lemme \ref{propzero}. Comme dans la d\'emonstration du th\'eor\`eme \ref{contractante} on obtient 
$$  \sum_{l=0}^{ N- [N\epsilon_{1}]}\varphi_{\alpha} (k-l) f(\frac{l}{N}) 
= O(N^{-\alpha-1}\epsilon_{1}^{-\alpha})=o(N^{-\alpha}).$$
On conclut en remarquant  que 
$$ \Bigl \vert \sum_{l= N- [N\epsilon_{1}]-1}^{N}\varphi_{\alpha} (N-l)  \left(f(\frac{l}{N})-f(1)\right)\Bigr \vert
=\le O( \frac{1}{N}) \sum (N-l)^{-\alpha}= 
O( N^{-\alpha}
\epsilon_{1}^{1-\alpha})= o(N^{-\alpha}).$$
\end{preuve}
  \subsection{Composition avec les d\'eriv\'ees enti\`eres}
 \begin{lemme}\label{lemmemoins1}
 Pour $N$ un entier naturel on consid\`ere $X_{N}$ le vecteur de $\mathbb R^{N+1}$
 d\'efini par $ X_{N}= \left( f(\frac{i}{N})\right)_{(i\in {0,\cdots, N})}$.  Pour tout entier $k$ v\'erifiant 
 $n<k<N-n$, et pour toute fonction  $f \in C^{n}]0,1[$ on a alors 
 $ \left(T_{N} (\varphi_{n}) X_{N}\right)_{k}=\frac{ (2)^{n}}{N^{n}} f^{(n)} (\frac{k}{n}) + o(\frac{1}{N^{n}})$.
 \end{lemme}
 \begin{preuve}{}
 On v\'erifie facilement que pour tout entier $q$ 
 on a pour $i,j$ v\'erifiant $q<i,j<N-q$
 \begin{equation}\label{compo}
 \left( T_{N}(\varphi_{1}) T_{N}(\varphi_{q-1})\right)_{i,j} 
 = T_{N}(\varphi_{q})_{i,j}.
 \end{equation}
 Si $n$ est maintenant un entier fix\'e, on a $n-1$ constantes $a_{i,1}, 2\le i \le n$ telles que 
  pour $1<j<N-1$,
 $$\left(T_{N}(\varphi_{1}) X_{N}\right)_{j}= (-2)\left( \frac{1}{N} f^{(1)} (\frac{j}{N}) 
 + \sum_{i=2}^{n} a_{i,1}
 f^{(i)} (\frac{j}{N})\right) + o (\frac{1}{N^{n}}).$$
 En utilisant l'\'equation (\ref{compo}) on obtient l'existence de $n-3$ constantes $a_{i,2},  3\le i \le n$ telles que 
pour $2<j'<N-2$ 
 $$
\left(T_{N}(\varphi_{2}) X_{N}\right)_{j'} =
  4 \frac{1}{N^{2}} f^{(2)} (\frac{j'}{N}) 
 + \sum_{i=3}^{n-1} a_{i,2}
 f^{(i)} (\frac{j'}{N})+  a_{n,1}\left(f^{(n)} (\frac{j'+1}{N}) - f^{(n)} (\frac{j'-1}{N})\right)+o (\frac{1}{N^{n}}),
$$
 ce qui donne 
 $$
 \left(T_{N}(\varphi_{2}) X_{N}\right)_{j'} =
  4 \frac{1}{N^{2}} f^{(2)} (\frac{j'}{N}) 
 + \sum_{i=3}^{n-1} a_{i,2}
 f^{(i)} (\frac{j'}{N})+ o (\frac{1}{N^{n}}).
$$
  En r\'ep\'etant le principe suffisamment de fois on obtient le lemme.
   \end{preuve}
 \begin{theoreme}\label{diff}
 Soit $\alpha \in ]0,1[$ et $n$ un entier. On suppose que $f$ est une fonction n fois d\'erivable   dans $[0,1]$ 
 v\'erifiant $f(0)=f^{(1)}(0)=\cdots= f^{(n-1)}(0)=0$,
 $ f(1)= f^{(1)}(1)= \cdots=f^{(n-1)}(1)=0$ et $f^{(n)}$ contractante.
 Alors  $f$ est $D^{\alpha+n}$ d\'erivable si et seulement si $f^{(n)}$ est localement contractante et dans 
 $L^1_{1-\alpha}([0,1])$. Alors pour tout r\'eel $x$ dans $]0,1[$ 
 $$ D^{\alpha+n} (f) (x) = 2^nD^{\alpha} (f^{(n)}) (x).$$
 \end{theoreme}
 \begin{remarque}
 Les hypoth\`eses sur $f^{(n)}$  permettent d'affirmer,
 gr\^{a}ce au th\'eor\`eme \ref{contractante} que $D_{\alpha} (f^{(n)}) (x)$ 
 est bien d\'efini pour $x \in ]0,1[$.
 \end{remarque}
 Cette propri\'et\'e et une cons\'equence directe du lemme ci-dessous
 \begin{lemme}
  Soit $\alpha \in ]0,1[$ et $n$ un entier. On suppose que $f$ d\'erivable
  dans $[0,1]$ 
 v\'erifiant $f(0)=f(1)=0$, et $f^{(1)}$ contractante.
 Alors  $f$ est $D^{\alpha+1}$ d\'erivable en $x\in ]0,1[$ si et seulement si $f^{(1)}$ est  $\alpha$ d\'erivable en $x$ et 
 $$ D^{\alpha+1} (f) (x) = (2) D^{\alpha} (f^{(1)}) (x).$$
 \end{lemme}
 \begin{preuve}{}
 Posons $k$ l'entier tel que $k=[Nx]$. On sait que 
$$ D^{\alpha+1} (f) (x) = 
N^{\alpha+1} \sum_{l=0}^N \widehat{\varphi_{\alpha+1}} (k-l)
f(\frac{l}{N}) +o(1).$$
Ce qui peut s'\'ecrire aussi 
$$ D^{\alpha+1} (f) (x) = 
N^{\alpha+1} \sum_{l=0}^N 
\left( 
\widehat{\varphi_{\alpha}} (k-l+1) \widehat{\varphi_{1}} (-1)+
\widehat {\varphi_{\alpha}} (k-l-1) \widehat{\varphi_{1}} (1)
\right)f(\frac{l}{N}) +o(1).$$
Ce qui s'\'ecrit encore
$$ D^{\alpha+1} (f) (x) = 
N^{\alpha+1} \sum_{h=1}^{N-1} \widehat{\varphi_{\alpha}} (k-h)
\left( \widehat{ \varphi_{1}}(-1) f(\frac{h+1}{N})+\widehat{ \varphi_{1}}(1) f(\frac{h-1}{N})\right)+
R_{N}+o(1)$$
avec $R_{N}= \widehat{\varphi_{\alpha}} (k-N) \widehat{\varphi_{1}}(1) f (\frac{N-1}{N})+
 \widehat{\varphi_{\alpha} }(k) \widehat{\varphi_{1}}(-1) f (\frac{1}{N})=o(N^{-\alpha-2}).$
 Puisque  $f^{(1)}$  est contractante nous pouvons \'ecrire, $c$ \'etant un r\'eel dans 
 $]\frac{h-1}{N},\frac{h+1}{N}[$  
 $$  \widehat{ \varphi_{1}}(-1) f(\frac{h+1}{N})+\widehat{\varphi_{1}}(1) f(\frac{h-1}{N})
 = \frac{2}{N} f^{(1)}(c)= \frac{2}{N} \left( f^{(1)}(\frac{h}{N}) + O(\frac{1}{N})\right)$$
 ce qui donne 
 \begin{align*} D^{\alpha+1} (f) (x) &= 
2N^{\alpha} \sum_{h=1}^{N-1} \widehat{\varphi_{\alpha}} (k-h)
\left( f^{(1)} (\frac{h}{N}) +O(\frac{1}{N}) \right) +o(1)\\
&= N^{\alpha} \sum_{h=0}^{N} \widehat{\varphi_{\alpha}} (k-h)
f^{(1)} (\frac{h}{N}) +o(1)
\end{align*}
avec $R'_{N}=-2 N^{\alpha}\left( \widehat{\varphi_{\alpha}} (k-N) f^{(1)}(1)+ 
\widehat{\varphi_{\alpha}} (k) f^{(1)}(0)\right)$.
Ce qui ach\`eve de prouver le lemme.
 La r\'eciproque s'obtient en remontant les calculs.
 \end{preuve}
 \begin{remarque}
 Si $\alpha \in ]0,1[$ et si  $f$ 
 est  une fonction d\'erivable de d\'eriv\'ee contractante sur $[0,1]$,  v\'erifiant de plus $f(0)=f(1)=0$, alors un calcul direct permet de v\'erifier que $D^{\alpha} (f)$ est d\'erivable sur $]0,1[$ et de d\'eriv\'ee $D^{\alpha} (f')$.
 Si on suppose de plus $f'(0)=f'(1) =0$ alors ce r\'esultat peut \^{e}tre \'etendu \`a 
 $[0,1]$ tout entier.
\end{remarque}

  \section{La formule d'inversion des matrices de Toeplitz \`a coefficients r\'eels}
  Classiquement, on sait calculer explicitement les coefficients de l'inverse  d'une matrice de Toeplitz de symbole positif $f$ qui v\'erifie $\ln f \in L^1 (\mathbb T)$. Ces formules 
  utilisent fortement le fait qu'une telle fonction $f$ peut s'\'ecrire $f = g \bar g$ avec 
  $g \in H^+ (\mathbb T)$ ( \cite{RRS}) o\`u
  $H^{+}(\mathbb T) = \{ h\in L^{2} (\mathbb T) \vert \hat h (s)=0 \forall s, s<0\}$. Certains probl\`emes, notamment la recherche des valeurs
  propres, n\'ecessitent d'inverser des matrices de Toeplitz dont le symbole 
  n'est pas n\'ecessairement positif (voir, par exemple, \cite{JMR01} \cite {RS020}).
  Pour ce faire, nous allons consid\'erer des fonctions $f$ qui ne s'annulent pas sur le tore  admettant une d\'ecomposition  $f=g_{1}g_{2}$, avec 
   $g_{1}, g_{1}^{-1} \in H^+(\mathbb T)$, $\bar g_{2}, {\bar g_{2}}^{-1}\in H^+(\mathbb T)$.
   Nous posons  $\Phi_{N} = \chi^{N+1} \frac{g_{1}}{g{2}}$ et 
   $\tilde \Phi_{N} = \chi^{-N-1} \frac{g_{2}}{g{1}}$, 
   o\`u $\chi$ d\'esigne la fonction d\'efinie par $ \theta \mapsto e^{i\theta}$.
    Nous  consid\'erons alors les op\'erateurs de Hankel, $H_{\Phi_{N}} $ et $H_{\tilde \Phi_{N}} $
   d\'efinis par 
   $$ H_{\Phi_{N}} : H^{+} (\mathbb T)\mapsto {H^{+} }^{\perp} (\mathbb T)\quad 
   \phi \mapsto \pi_{-} ( \Phi_{N} (\phi) $$
   $$ H_{\tilde \Phi_{N}} : {H^{+} }^{\perp}(\mathbb T) \mapsto H^{+} (\mathbb T)
   \quad 
   \psi\mapsto \pi_{+} (\tilde  \Phi_{N}) (\psi). $$ 
   
   Dans la suite nous d\'esignerons par $\mathcal P_{N}$ l'ensemble des polyn\^{o}mes 
   trigonom\'etriques de degr\'e inf\'erieur ou \'egal \`a $N$, et nous notons 
   $\pi_{+}$ la projection orthogonale de $L^{2}(\mathbb T)$ dans $H^+ (\mathbb T)$ 
   et  $\pi_{+}$ la projection orthogonale de $L^{2}(\mathbb T)$ dans ${H^{+} }^{\perp}$.
   Enfin si $g\in L^{2}(\mathbb T)$ nous notons par 
   $d_{\infty}(g,\mathcal P_{N})$ la quantit\'e 
   $\inf_{q\in \mathcal P_{N}} \Vert g-q\Vert_{\infty}$, o\`u 
   $ \Vert g-q\Vert_{\infty} = \sup_{\theta\in [0, 2\pi[} \vert g (\theta) -q (\theta)\vert$.
  Nous pouvons alors \'enoncer le th\'eor\`eme :
  \begin{theoreme}\label{THINV}
  Si la fonction $f$ est comme ci-dessus, on suppose que les fonctions $g_{1}$ et $g_{2}$ v\'erifient 
  $$ 
  \lim_{N\rightarrow + \infty} d_{\infty}
  \left( \frac{1}{\vert g_{2}\vert ^{2}}, \mathcal P_{N}\right) =0
  \quad 
  \mathrm{et}
  \quad 
  d_{\infty}
  \left( \frac{1}{\vert g_{1}\vert ^{2}}, \mathcal P_{N}\right) < + \infty
  $$
  ou
   $$ 
  \lim_{N\rightarrow + \infty} d_{\infty}
  \left( \frac{1}{\vert g_{1}\vert ^{2}}, \mathcal P_{N}\right) =0
  \quad 
  \mathrm{et}
  \quad 
  d_{\infty}
  \left( \frac{1}{\vert g_{2}\vert ^{2}}, \mathcal P_{N}\right) <+ \infty
$$
  alors, pour $N$ suffisamment grand,
  \begin{enumerate}
  \item
   l'in\'egalit\'e $\Vert H_{\tilde \Phi_{N}} H_{\Phi_{N}} \Vert_{2}<1$ est v\'erifi\'ee (ce qui assure bien s\^{u}r l'inversibilit\'e de l'op\'erateur 
  $I - H_{\tilde \Phi_{N}} H_{\Phi_{N}} $),
  \item
  la matrice $T_{N}(f)$ est inversible,
  \item
    pour tout polyn\^{o}me $Q$ dans $\mathcal P_{N}$, on a 
 $$
  T_{N}^{-1} (f) (Q) =
  \pi_{+}(Q g_{2}^{-1}) g_{1}^{-1} - \pi_{+}
  \left ( \left(( I- H_{\tilde \Phi_{N}} H_{\Phi_{N}} )^{-1} \pi_{+} 
  \left( \pi_{+} (Q g_{2}^{-1}) \tilde \Phi_{N}\right)\right) \Phi_{N}\right ) g_{1}^{-1}.$$
  \end{enumerate}
  \end{theoreme} 
  Ce th\'eor\`eme a \'et\'e obtenu dans \cite{RRQuebl} pour inverser des matrices \`a blocs que nous interpr\'etons comme des matrices de Toeplitz tronqu\'ees dont le symbole 
  est une fonction matricielle contenue dans $L^2_{\mathcal M_{n}} (\mathbb T)
  =\{ M : \mathbb T \rightarrow \mathcal M_{n} (\mathbb C) \Big\vert 
 \int_{\mathbb T}\Vert M\Vert_{2}^2 d\sigma< + \infty\},$ o\`u 
 $\Vert M \Vert_{2}= [ \trace (M M ^*]^{1/2}$ et o\`u $\sigma$ est la mesure de Lebesgue su le tore $\mathbb T$. N\'eanmoins, il est facile de v\'erifier 
  qu'avec les hypoth\`eses que nous nous sommes donn\'ees (voir \cite{RRQuebl}, lemme 7.1.5 et corollaire 7.1.6), la d\'emonstration se transpose
  sans difficult\'e au cas des matrices de Toeplitz dont le symbole est une fonction de 
  $\mathbb T$ dans $\mathbb C$. \\
  Nous posons ici :
 $g_{1,R} = (1+R \chi)^\alpha$ et 
 $g_{2,R}= (1-R \bar \chi)^{\alpha}$, 
 et nous allons inverser la matrice 
 $T_{N,\alpha,R}$ de symbole $g_{1,R} g_{2,R} =
  (1-R^2 +R \chi -R\chi)^{2\alpha}$.
  Nous noterons $T_{N,\alpha,R}$ la matrice 
  $T_{N}(\phi_{\alpha,R})$.
\section{Calcul des coefficients de $\left(T^{-1}_{N,\alpha,R}\right)$}
\subsection{Notations}
Dans la suite de ce paragraphe nous supposerons $1>\alpha>0.$
On note  $\beta_{u,1,R}$ et $\beta_{u,2,R}$ les coefficients de Fourier 
de $g_{1,R} ^{-1}$ et $g_{2,R}^{-1}$. Nous avons \'evidemment pour $u\ge 0$ 
$ \beta_{u,1,R} = (-1)^{u} R^{u} \beta_{u}^{(\alpha)}$ et 
$ \beta_{u,2,R} =  R^{u} \beta_{-u}^{(\alpha)}$ avec (voir \cite{Zyg})
\begin{equation}\label{eq0}
 \beta_{u}^{(\alpha)} = \frac{u^{\alpha-1}}{\Gamma (\alpha)} 
+\frac{u^{\alpha-2}}{\Gamma (\alpha-1)}+o(u^{\alpha-2}).
\end{equation}
D'apr\`es le th\'eor\`eme \ref{THINV}
nous avons pour $0\le k,l\le N$ 
$$ \left(T^{-1}_{N,\alpha,R}\right)_{k+1,l+1} = T_{1,N,k,l,R}+ T_{2,N,k,l,R}$$
avec 
$$ T_{1,N,k,l,R}= \Bigl\langle \pi_{+} \left( \frac{\chi^{l}}{g_{2,R}}\right)\vert 
\pi_{+} \left( \frac{\chi^{k}}{\overline{g_{1,R}}}\right)\Bigr \rangle
$$
et 
$$ T_{2,N,k,l,R}= \Big\langle \pi_{+}
  \left(( I- H_{\tilde \Phi_{N,R}} H_{\Phi_{N,R}} )^{-1} \pi_{+} 
  \left( \pi_{+} (\frac{\chi^{l}}{ g_{2,R}}) \tilde \Phi_{N,R}\right)\right)\vert 
  \pi_{+} \left(\overline{ \Phi_{N}}\pi_{+}\left(\frac{\chi^k} {\overline {g_{1,R}}}\right)\right)
  \Bigr \rangle.
  $$
  Pour all\'eger les notations  nous noterons dans la suite du texte  
  $T_{1,k,l,R}$ pour $T_{1,N,k,l,R}$ et 
  $T_{2,k,l,R}$ pour $T_{2,N,k,l,R}$.
\subsection{Calcul de $T_{1,k,l,R}$}

Supposons tout d'abord : $k= [Nx]$, $l=[Ny]$, avec $0<x<y\le1$.\\
Avec les notations ci-dessus nous avons 
\begin{align*}
T_{1,k,l,R} &= \Bigl \langle \pi_{+} \left(\sum_{u\ge 0} \beta_{u}^{(\alpha)} (-1)^{u}R^{u} \chi^{l-u}\right)
\vert\pi_{+} \left( \sum_{v\ge 0}  R^{v} \overline{\beta_{v}^{(\alpha)} } \chi^{k-v}
\right) \Bigr \rangle\\
&= \Bigl \langle \sum_{u= 0}^{l} (-1)^{u}\beta_{u}^{(\alpha)} R^{u} \chi^{l-u}
\vert \sum_{v= 0} ^{k} R^{v} \overline{\beta_{v}^{(\alpha)} } \chi^{k-v}\Bigr \rangle.
\end{align*}
Ce qui donne 
\begin{align*}
 T_{1,k,l,R} &= \sum_{u=0}^{k} R^{l-k+2u} (-1)^{l-k+u} \beta_{l-k+u}^{(\alpha)} 
 \overline{\beta_{u}^{(\alpha)}}\\
 &=(-R)^{l-k} \left(\sum_{u=0}^{k_{0}} R^{2u} (-1)^{l-k+u}  \beta_{l-k+u}^{(\alpha)} +\sum_{u=k_{0}+1}^{k} R^{2u} (-1)^{u}  \beta_{l-k+u}^{(\alpha)} 
 \overline{\beta_{u}^{(\alpha)}}\right)
 \end{align*}
o\`u $k_{0}=N^{\beta}$ o\`u $\beta$ est un r\'eel strictement compris entre $0$ 
  et $1$. 
  Une sommation d'Abel permet de v\'erifier que 
  $$ \sum_{u=k_{0}+1}^{k} R^{2u} (-1)^u \beta_{l-k+u}^{(\alpha)} 
 \overline{\beta_{u}^{(\alpha)}} = O\left( (l-k)^{\alpha-1} 
 k_{0}^{\alpha-1} \right) = O\left( N^{(\beta+1)(\alpha-1)}\right).
 $$
   En posant  
   $ \psi_{1} (l-k,u) =  \beta_{l-k+u}^{(\alpha)}- \frac{(l-k)^{\alpha-1}} {\Gamma (\alpha)} $ 
  on peut \'ecrire 
  \begin{align*}
   \sum_{u=0}^{k_{0}} R^{2u} (-1)^u \beta_{l-k+u}^{(\alpha)} 
 \overline{\beta_{u}^{(\alpha)}} &= \frac{(l-k)^{\alpha-1} }{\Gamma (\alpha)} 
 \sum_{u=0}^{k_{0}} R^{2u} (-1)^u \overline{\beta_{u}^{(\alpha)}} 
 \psi_{1} \left( l-k,u\right)\\
 &+ \sum_{u=0}^{k_{0}} (-1)^u R^{2u} \overline{\beta_{u}^{(\alpha)}}  \frac{(l-k)^{\alpha-1}}{\Gamma (\alpha)} 
 \end{align*}
 ce qui donne 
 $$ \sum_{u=0}^{k_{0}} R^{l-k+2u} (-1)^u \beta_{l-k+u}^{(\alpha)} 
 \overline{\beta_{u}^{(\alpha)}} = R^{l-k} \frac{(l-k)^{\alpha-1} }{\Gamma (\alpha)} 
 \sum_{u=0}^{k_{0}} R^{2u} (-1)^u \overline{\beta_{u}^{(\alpha)}} 
 + O \left( (l-k)^{\alpha-2} k_{0}^\alpha \right).
 $$
 Enfin puisque
  $\displaystyle{\sum_{u>k_{0}} R^{2u}(-1)^{u} \beta_{u}^{(\alpha)}}=
   O(k_{0}^{\alpha-1})$ nous pouvons \'enoncer la propri\'et\'e suivante :
 \begin{prop}\label{propo1}
En supposant $k= [Nx]$, $l=[Ny]$, avec $0< x\neq y\le1$ on a, pour tout r\'eel $\beta \in ]0,1[$ et pour $N$ assez grand : 
 $$ T_{1,k,l,R} = (-1)^{l-k}  R^{\vert l-k\vert } \frac{\vert l-k\vert ^{\alpha-1} }{\Gamma (\alpha)} 
 (1+R^2)^{\alpha} \left(1+R_{1}\vert l-k\vert \right) \quad \mathrm{avec} \quad 
 \vert R_{1}\vert l-k\vert \vert \le M_{1} (l-k)^{ \beta(\alpha-1)},$$
 $M_{1}$ \'etant une constante ind\'ependante de $k$ et $l$.
 \end{prop}
  \subsection{Estimation de $T_{2,k,l,R}$}
  On a (voir l'appendice pour la d\'emonstration
  de ce r\'esultat)
\begin{prop} \label{propo3}
Pour tout intervalle $[\delta_{1},\delta_{2}]$ contenu dans $[0,1]$ 
$T_{2,k,l,R} =o(N^{\alpha-1})$ uniform\'ement pour 
tout $(x,y)$ dans $[\delta_{1}, \delta_{2}]\times [0,1]$,
\end{prop}
\section{Inversion de l'op\'erateur diff\'erentiel $D_{\alpha}$
pour $\alpha\in ]0,1[$.}
  \subsection{Deux lemmes techniques}
  Pour la suite de nos calcul nous allons devoir utiliser 
  ls deux lemmes suivants. 
  \begin{lemme} \label{LEMATA}
  Soit $\beta$ un r\'eel dans $]0,1[$.
  
  Si $\vert l-k\vert > N^{\beta}$
  alors pour tout $\beta_{1}\in ]0,1[$ et pour $N$ assez grand 
  on a 
  $$ T_{1,k,l,R}= (-1)^{\vert l-k\vert } R^{\vert l-k\vert } 
  \frac{\vert l-k\vert^{\alpha-1}} {\Gamma (\alpha)}
  (1+R^2)^{-\alpha} \left(1+R_{4}(\vert l-k\vert )\right)$$
  avec $\vert R_{4}(l-k)\vert <  K_{4 } (l-k)^{\beta_{1}(\alpha-1)}$ 
  le r\'eel  
  $K_{4}$ \'etant une constante ind\'ependante de $k$, $l$ 
  et $R$.
   \end{lemme}
  Ce lemme est en fait une variante de la  propre\'et\'e \ref{propo1} et il se d\'emontre de la m\^{e}me mani\`ere.
    \begin{lemme}\label{LEMATA2}
  Soit $\alpha \in ]0,1[$ et $f $ une fonction localement contractante sur $]0,1[$ et dans $L^{1} _{1-\alpha}([0,1])$ 
  pour $0\le \gamma\le 1-\gamma$. 
  Alors pour un entier $N$ fix\'e suffisamment grand, et pour $k =[Nx]$, $0<x<1$
  $$\lim_{R\rightarrow 1} \sum_{l=0}^N T_{1,k,l,R} f(\frac{l}{N})
  = N^{\alpha} \frac{1}{2^\alpha \Gamma (\alpha)}
   \int_{0}^x \frac{f(t)}{(t-x)^{1-\alpha} }dt+ o(N^\alpha).$$
  \end{lemme}
  \begin{preuve}{}
  On peut \'ecrire 
  \begin{align*}
  \sum_{l=0}^N  T_{1,k,l,R} f(\frac{l}{N}) &= 
  \sum_{l=0}^{k-k_{0}-1}  T_{1,k,l,R} f(\frac{l}{N})+
  \sum_{l=k-k_{0}}^{k+k_{0}} T_{1,k,l,R} f(\frac{l}{N})\\
  &+
  \sum_{l=k+k_{0}+1}^{N} T_{1,k,l,R} f(\frac{l}{N})
  \end{align*}
  o\`u $k_{0}= N^\beta$, avec $\beta\in ]0,1[$ (la valeur de $\beta$ sera pr\'ecis\'ee ult\'erieurement). 
  En utilisant le calcul de $T_{1,k,l,R}$ et le lemme \ref{LEMATA} 
  on obtient si $ \beta_{1}$ comme dans le lemme 
  \ref{LEMATA} 
  $$ \sum_{l=k+k_{0}+1}^{N}  T_{1,k,l,R} f(\frac{l}{N}) 
  = \sum_{l=k+k_{0}+1}^{N} \left( (-1)^{k-l} R^{k-l} \frac{\vert k-l\vert ^{\alpha-1}} {\Gamma (\alpha)}
  (1+R^2)^{-\alpha} f(\frac{l}{N}) \right)+ 
  O \left( (l-k)^{(\alpha-1)(\beta_{1}+1)}\right)
$$
ce qui donne, si $\beta_{1}$ bien choisi,
$$ \sum_{l=k+k_{0}+1}^{N}  T_{1,k,l,R} f(\frac{l}{N}) 
  = \sum_{l=k+k_{0}+1}^{N} \left( (-1)^{k-l} R^{k-l} \frac{\vert k-l\vert ^{\alpha-1}} {\Gamma (\alpha)}
  (1+R^2)^{-\alpha} f(\frac{l}{N}) \right)+ 
 o(N^{\alpha})
$$
uniform\'ement par rapport \`a $R$. D'o\`u :
$$ \lim_{R \rightarrow 1} \sum_{l=k+k_{0}+1}^{N}  T_{1,k,l,R} f(\frac{l}{N}) 
= \frac{1} {2^{\alpha} \Gamma \alpha}\sum_{l=0}^{k-k_{0}-1} \left( (-1)^{k-l}  \vert k-l\vert ^{\alpha-1} f(\frac{l}{N}) \right)
 +o(N^{-\alpha}).
$$
Une sommation d'Abel permet d'\'ecrire, puisque on a, par convention, pos\'e $f(N) =f(0) =0$,
\begin{align*}
 \Bigl \vert\sum_{l=k+k_{0}+1}^{N-1} (-1)^{k-l}  \vert k-l\vert^{\alpha-1} f(\frac{l}{N}) 
\Bigr \vert 
&\le \Bigl \vert(N-k)^{\alpha-1} f(1) +(k_{0}+1)^{\alpha-1} f(\frac{k+k_{0}+1}{N})\Bigr\vert\\
&+ 
\sum_{l=k+k_{0}+1}^{N-1} 
\Bigl\vert \vert k-l+1\vert ^{\alpha-1} f(\frac{l+1}{N}) - \vert k-l\vert^{\alpha-1} f(\frac{l}{N})
\Bigr \vert.
\end{align*}
On \'ecrit 
\begin{align*}
&\sum_{l=k+k_{0}+1}^{N-1} 
\Bigl\vert \vert k-l+1\vert ^{\alpha-1} f(\frac{l+1}{N}) - \vert k-l\vert ^{\alpha-1} f(\frac{l}{N})\Bigr\vert \le
\sum_{l=k+k_{0}+1}^{N-1} 
\Bigl\vert \vert k-l+1\vert ^{\alpha-1} f(\frac{l+1}{N}) 
\\ &-\vert k-l\vert ^{\alpha-1} f(\frac{l+1}{N})\Bigr\vert 
+\Bigl\vert \vert k-l\vert ^{\alpha-1} f(\frac{l+1}{N}) - \vert k-l\vert ^{\alpha-1} f(\frac{l}{N})\Bigl\vert.
\end{align*}
Les hypoth\`eses faites sur $f$ donnent, si $k' =[Nx']$ avec $x<x'<1$.
\begin{align*}
\sum_{l=[Nx']}^{N-1} 
\Bigl\vert \vert k-l+1\vert ^{\alpha-1} f(\frac{l+1}{N}) 
\\ &-\vert k-l\vert ^{\alpha-1} f(\frac{l+1}{N})\Bigr\vert \\
&\le O(N^{\gamma}) \sum_{l=[Nx']}^{N-1} \vert k-l\vert ^{\alpha-1} = O(N^{\gamma+\alpha-1})
=o(N^{\alpha}).
\end{align*}
Et d'autre part 
$$ \sum_{l=k+k_{0}+1}^{[Nx']-1} 
\Bigl\vert \vert k-l+1\vert ^{\alpha-1} f(\frac{l+1}{N}) 
-\vert k-l\vert ^{\alpha-1} f(\frac{l+1}{N})\Bigr\vert
=o(N^{\alpha}).$$
Enfin la somme
$$ \Bigl\vert \vert k-l\vert ^{\alpha-1} f(\frac{l+1}{N}) - \vert k-l\vert ^{\alpha-1} f(\frac{l}{N})\Bigl\vert$$ 
a \'et\'e trait\'ee dans la d\'emonstration du th\'eor\`eme \ref{contractante}, on sait qu'elle est d'ordre $o(N^{\alpha})$.
Ce qui implique finalement 
$$
\sum_{l=k+k_{0}+1}^{N}
\vert (k-l+1)^{\alpha-1} f(\frac{l+1}{N}) - (k-l)^{\alpha-1} f(\frac{l}{N})\vert
= o(N^{\alpha}).
$$
Et on peut conclure 
$$ \lim_{R\rightarrow 1}\sum_{l=0}^{k-k_{0}-1}  T_{1,k,l,R} f(\frac{l}{N}) =o(N^{\alpha}).$$
On a d'autre part, si $\beta$ choisi de tel mani\`ere que
$ 0<\beta<\alpha$ :
$$ \sum_{l=k-k_{0}}^{k+k_{0}}  T_{1,k,l,R} f(\frac{l}{N}) =O(N^{\beta})=o(N^{\alpha}).$$
Et enfin, comme pour la pr\'ec\'edente somme  
$$ \sum_{l=0}^{k-k_{0}-1} T_{1,k,l,R} f(\frac{l}{N})
  = \sum_ {l=0}^{k-k_{0}-1}\left( R^{k-l} \frac{(k-l)^{\alpha-1}} {\Gamma (\alpha)}
  (1+R^2)^{-\alpha} f(\frac{l}{N}) \right)+o(N^{\alpha}),
$$
uniform\'ement par rapport \`a $R$,
soit
$$ \lim_{R\rightarrow 1} \sum_{l=0}^{l=k-k_{0}-1}T_{1,k,l,R} f(\frac{l}{N})
=  \frac{1}
  {\Gamma (\alpha) 2^{\alpha}}\sum_{l=0}^ {l=k-k_{0}-1}
  \left(  (l-k)^{\alpha-1} f(\frac{l}{N}) \right)+O( N^{\alpha}).$$
  Ce qui donne facilement 
\begin{equation}\label{LIM} 
\lim_{R\rightarrow 1} \sum_{l=0}^{l=k-k_{0}-1} T_{1,k,l,R} f(\frac{l}{N})
= N^{\alpha}\frac{1} {\Gamma (\alpha)2^{\alpha}} \int_{0}^{x} (t-x)^{\alpha-1}f(t) dt 
+o(N^{\alpha}).
\end{equation}
  \end{preuve}
  \subsection {D\'efinition de l'op\'erateur $J_{\alpha}$ pour $0<\alpha<1$.}
  \begin{definition} 
 Soit $f$ une fonction d\'efinie sur $]0,1[.$ Pour  $ x \in ]0,1[ $ et $\alpha$ un r\'eel positif  on notera 
 $J_{\alpha} \left( f(x)\right)$ le r\'eel, s'il existe 
 $$ \lim_{N\rightarrow + \infty, R\rightarrow 1^{-}} 
 N^{-\alpha} \sum_{l=0}^{N} 
 \left(T_{N,\alpha,R}^{-1}\right)_{k+1,l+1} f\left( \frac{l}{N}\right)
 =J_{\alpha} (x)$$
 avec $\lim_{N\rightarrow + \infty} \frac{k}{N}=x$.
\end{definition}
En utilisant le lemme \ref{LEMATA2} on peut \'ecrire 
\begin{prop} \label{INTFRAC}
Si $\alpha \in ]0,1[$ et $f $ une fonction localement contractante sur $]0,1[$ et dans $L^{1}_{1-\alpha} (]0,1[)$ alors pour $x\in ]0,1[ $ on a 
  $$ J_{\alpha} (x) = \frac{1}{2^\alpha \Gamma (\alpha)}
  \int_{0}^x \frac{ f(t)} {(t-x)^{1-\alpha}} dt.$$
\end{prop}
  \subsection{R\'esolution d'une \'equation diff\'erentielle fractionnaire dans $]0,1[$ pour 
  $\alpha\in]0,1[$.}
   \begin{prop}\label{invgauche} 
   Soient une fonctions $f$  localement contractantes sur $]0,1[$, appartenant  \`a $L_{1-\alpha}^1([0,1])$  et  telle
   qu'il existe une fonction $\psi$ d\'efinie et born\'ee sur $]0,1[$ telle que  que  $\psi = D_{\alpha} (f)$ sur $[0,1]$.    Alors pour tout  $x\in ]0,1[$ la quantit\'e $\left( J_{\alpha} (\psi) \right)(x)$ 
   est d\'efinie et on a 
  $f(x)=\left( J_{\alpha} (\psi) \right)(x)$. 
  \end{prop}
 
  \begin{preuve}{}
  
 Soit $N$ fix\'e. D'apr\`es la d\'emonstration du th\'eor\`eme
  \ref{contractante} en choisissant un r\'eel $\epsilon_{1}>0 $ assez petit
 on a , uniform\'ement pour  $t$ v\'erifiant 
  $1-\epsilon_{1}>t>\epsilon_{1}$ 
  \begin{equation} \label{ronron}
  N^\alpha \sum_{j=0}^N \left(T_{N}
  ( \varphi_{\alpha})\right)_{k+1,j+1} f(\frac{j}{N}) -
  \left(D_{\alpha} f\right)(t) = o(1)
  \end{equation}
 D\'efinissons maintenant les vecteurs $X_{N}$ et 
 $Y_{N}$ de $\mathbb R^{N+1}$ de la mani\`ere suivante.
 $$ (X_{N})_{(0)}=(X_{N})_{(1)}=0, (X_{N})_{k} = f(\frac{k}{N}), 
 \quad \rm{pour} \quad 1\le k \le N-1, $$
 $$ (Y_{N})_{(0)}=(Y_{N})_{(1)}=0, (Y_{N})_{k} = \psi(\frac{k}{N}), 
 \quad \rm{pour} \quad [N\epsilon_{1}]\le k \le N-[N\epsilon_{1}], $$
 et enfin 
 $(Y_{N})_{k}= N^{\alpha} \displaystyle{\sum_{l=0}^N 
 \left(T_{N,\alpha,R}\right)_{k+1,l+1} (X_{N})_{l}}$
 pour $0\le l\le [N\epsilon]-1$ et
 $ N-[N\epsilon]+1\le l \le N$.
 Nous pouvons \'ecrire, en notant $ T_{N,\alpha,R}$ la matrice 
 $T_{N}(\varphi_{\alpha,R} )$,
 $$ N^\alpha \left( T_{N,\alpha,R} (X_{N}) \right)
 = Y_{N}+R_{N}$$
 avec, en utilisant l'\'equation (\ref{ronron}),
 $(R_{N})_{k}= o(1)$ uniform\'ement pour $[N\epsilon_{1}]\le k \le N-[N\epsilon-1]$ et $(R_{N})_{k} = 0$ pour $0\le k\le 
 [N\epsilon_{1}]-1$
 et $ N-[N\epsilon]+1\le k\le N$.
 On a \'evidemment pour $k$ compris entre $[N\epsilon_{1}]$
 et $N-[N\epsilon_{1}]$,
  $$ (X_{N})_{k}=
 N^{-\alpha}  \left(T_{N\alpha,R} ^{-1} ( Y_{N}+ R_{N})\right)_{k}.$$
Soit maintenant $x\in ]0,1[$. On peut se ramener \`a $x\in [\epsilon_{1}, 1-\epsilon_{1}]$. On a  $[Nx]=O(N)$ ce qui implique, pour $N$ assez grand :
$N-[N\epsilon_{1}]>[Nx]>[N\epsilon_{1}].$
  Si $k=[Nx]$ nous avons donc :
  $$ (X_{N})_{k}= N^{-\alpha}
  \left( T_{N,\alpha,R}^{-1}Y_{N}\right)_{k}
  + \left(  T_{N,\alpha,R}^{-1}R_{N}\right)_{k}.$$
  Nous pouvons \'ecrire 
  \begin{align*}
  N^{-\alpha}  \left( \left( T_{N} (\phi_{\alpha}) \right)^{-1}Y_{N}\right)_{k} 
  &=   N^{-\alpha}  \sum_{l=0}^{[N\epsilon]-1} \left(T_{N,\alpha,R} ^{-1}\right)_{(k+1,l+1)} (Y_{N})_{l}\\
  &+   N^{-\alpha} \sum_{[N\epsilon]}^{N-[N\epsilon]}
   \left(T_{N,\alpha,R} ^{-1}\right)_{(k+1,l+1)} (Y_{N})_{l}\\
 & +   N^{-\alpha} \sum_{N-[N\epsilon]+1}^N
 \left(  T_{N,\alpha,R}^{-1}\right) _{(k+1,l+1)} (Y_{N})_{l}.
   \end{align*}
  L'\'etude des coefficients de $  T_{N,\alpha,R}^{-1}$ permet d'\'ecrire,   $$   N^{-\alpha}  \sum_{l=0}^{[N\epsilon]-1}
   \left(T_{N,\alpha,R} ^{-1}\right)_{(k+1,l+1} (Y_{N})_{l} =O(\epsilon_{1})$$
  et 
 $$ N^{-\alpha} \sum_{l=N-[N\epsilon]+1}^N
 \left(T_{N,\alpha,R}^{-1}\right)_{(k+1,l+1} (Y_{N})_{l}
  =O(\epsilon_{1}).$$
 D'autre part on a aussi, toujours gr\^{a}ce \`a l'\'etude
 des coefficients de $ \left( T_{N,\alpha,R} \right)^{-1}$ :
  $$ N^{-\alpha} \sum_{l=[N\epsilon]}^{N-[N\epsilon]}
  \left(T_{N,\alpha,R}^{-1}\right)_{(k+1,l+1} (R_{N})_{l}
  =o(1).$$
  C'est \`a dire que 
  $$ (X_{N})_{k}= N^{-\alpha} \sum_{l=[N\epsilon]}^{N-[N\epsilon]}
 \left(T_{N,\alpha,R} ^{-1}\right)_{(k+1,l+1} (Y_{N})_{l}
  +o(1)+O(\epsilon_{1}).$$
  Par des arguments identiques \`a ceux d\'ej\`a utilis\'es et en consid\'erant $\epsilon_{1}$ suffiasamment petit on obtient 
  $$ (X_{N})_{k}= N^{-\alpha} \sum_{l=0}^{N}
   \left( T_{N,\alpha,R}^{-1}  \right)_{(k+1,l+1} (Y_{N})_{l}
  + o(1).$$
  Ce qui implique 
  $\psi (x) = \left(J_{\alpha}(f)\right)(x),$  en passant \`a la limite sur $N$ et $R$.
  \end{preuve}
    Nous obtenons imm\'ediatement le corollaire suivant : 
  \begin{corollaire}\label{equdiff1}
  Soit $\alpha\in ]0,1[$ et $\psi$ une fonction localement contractante sur $]0,1[$ 
  appartenant \`a $L^1_{1-\alpha} ([a,b])$. Alors la fonction $J_{\alpha} (f)$ d\'efinie sur $[0,1]$ par 
  $ \left(J_{\alpha}(\psi)\right)(x)= \frac{1}{2^{\alpha}\Gamma\alpha)}
  \int_{0}^{x} \frac{\psi(t)}{(x-t)^{1-\alpha}} dt$ est d\'efinie sur $[0,1]$, localement contractante sur $]0,1[$, appartient \`a l' ensemble $L^1_{1-\alpha}([a,b])$, et est solution sur $]0,1[$ de l'\'equation diff\'erentielle $D_{\alpha} (y) =\psi$.
  \end{corollaire}
  \begin{preuve}{} 
  Il est clair que la fonction $J_{\alpha}(\psi)$ est d\'efinie sur $[0,1]$.
  Cette fonction est \'egalement localement contractante sur $[0,1]$. En effet 
   notons par $K$ le r\'eel positif tel que 
  pour tout $x,x'\in [\delta_{1},\delta_{2}]\subset ]0,1[$ on ait $\vert \psi(x') -\psi(x)\vert \le K \vert x'-x\vert$.
  En effet, si on suppose $\delta_{1}\le x\le x'\le \delta_{2}$ on a 
  $$ \vert \left( J_{\alpha}(\psi)\right) (x') - \left( J_{\alpha}(\psi)\right) (x)\vert 
  \le \Bigl\vert \int_{0}^{x} \left(\psi (x'-u) -\psi(x-u) \right) u^{\alpha-1}du \Bigr\vert 
  +\Bigl \vert \int_{x}^{x'} \psi (x'-u) u^{\alpha-1} du\Bigr \vert.$$
  La valeur absolue de la premi\`ere int\'egrale est major\'ee par 
  $ K \vert x'-x\vert  \int_{0}^{1} u^{\alpha-1} du $.
  la seconde v\'erifie 
  \begin{align*}
  \vert \int_{x}^{x'} \psi (x'-u) u^{\alpha-1} du\vert&\le
   \int_{x}^{x'} \vert \psi(x'-u) -\psi(0)\vert u^{\alpha-1} du \\
   &\le  K \int_{x}^{x'} (x'-u) u^{\alpha-1} du \le K (x'-x)\int_{0}^{1} u^{\alpha-1} du 
   \end{align*}
  Ceci suffit \`a prouver que la fonction $J_{\alpha}(\psi)$ est  contractante sur 
  $[\delta_{1},\delta_{2}]$. On v\'erifie \'egalement facilement que 
 $J_{\alpha} (\psi)\in L^{1}_{1-\alpha}([0,1])$. 
   La quantit\'e $\left(D_{\alpha}\left( J_{\alpha}(\psi)\right) \right) (x)$ est donc bien d\'efinie 
   pour tout $x$ dans $]0,1[$. Si $\psi' = D_{\alpha}\left( J_{\alpha}(\psi)\right)$, 
   En utilisant la propri\'et\'e  \ref{INTFRAC} on sait  que $J_{\alpha}(\psi')$ est d\'efini et vaut $J_{\alpha}(\psi)$, autrement dit 
   $J_{\alpha}(\psi-\psi')=0$. Cette fonction est donc ind\'efiniment d\'erivable sur 
   $]0,1[$, ce qui ne se peut que si $\psi=\psi'$. 
    \end{preuve}
\begin{remarque}
Puisque la fonction nulle est une fonction contractante 
sur $[0,1]$ qui vaut z\'ero en z\'ero on d\'eduit du corollaire 
\ref{equdiff1} l'\'equivalence $ D_{\alpha}(f) (x) =0 \iff f(x)=0
\forall x \in ]0,1[$.
\end{remarque}
       \subsection{ L'op\'erateur $J_{n}$ pour $n$ entier}
 \begin{definition} Pour toute fonction $f$ dans $C([0,1)]$, et tout entier naturel non nul pair $n=2p$ on pose  
  $ J_{n} (f) = (-1)^p \int_{0}^1G_{p}(x,t) f(t) dt$,
  avec , pour $0\le x,y\le 1$ 
  $$G_{p}(x,y) = \frac{x^py^p} {\left((p-1)!\right)^2} 
  \int_{\max (x,y)}^1\frac{(t-x)^{p-1} (t-y)^{p-1}}{t^{2p}} dt
  \quad \mathrm{si} \quad (x,y)\neq (0,0)$$
  et $G(0,0)=0$.
  \end{definition}
  On rappelle qu'il a \'et\'e \'etabli dans (\cite{RS04}) que 
  $G_{p}$ est le noyau de Green associ\'e \`a l'op\'erateur 
  diff\'erentielle $ (-1)^p d^{2p}/dx^{2p} (f) $ sur $[0,1]$ associ\'e 
  aux conditions initiales 
  $$ f^{(0)} (0) = f^{(1)} (0)= \cdots = f^{(p-1)} (0)=0,
  \quad 
   f^{(0)} (1) = f^{(1)} (1)= \cdots = f^{(p-1)} (1)=0.$$
    Nous pouvons alors \'enoncer les lemmes suivants 
   \begin{lemme}
   Si $n$ est un entier pair non nul et si $f$ est une fonction continue sur $[0,1]$ alors 
   $$ (-1)^{n/2} D_{n}\circ J_{n}(f) = f$$
  \end{lemme}
  C'est une cons\'equence imm\'ediate de la d\'efinition de 
  $J_{n}$ et du rappel ci-dessus.
  \begin{lemme}  \label{previous2} 
   Si $n$ est un entier pair non nul et si $f$ est une fonction de classe $C^n$ sur $[0,1]$v\'erifiant les conditions initiales 
   $$ f^{(0)} (0) = f^{(1)} (0)= \cdots = f^{(p-1)} (0)=0,
  \quad 
   f^{(0)} (1) = f^{(1)} (1)= \cdots = f^{(p-1)} (1)=0,$$ 
   on a 
      $$ J_{n} \circ D_{n}(f)=f.$$
  \end{lemme}
  \begin{preuve}{}
  D'apr\`es le lemme pr\'ec\'edente 
  $ D_{n} \left( J_{n}\circ D_{n} (f) \right) =
  D_{n}(f)$. Les fonctions $J_{n}\circ D_{n} (f)$ et 
  $f$ diff\`erent donc sur $[0,1]$ d'un polyn\^{o}me de degr\'e $n-1$.
  Mais d'apr\`es la d\'efinition de $J_{n}$ et les conditions 
  impos\'ees \`a $f$ ce polyn\^{o}me est n\'ecessairement nul.
  \end{preuve}
  \subsection{L'op\'erateur $\tilde J_{\alpha}$}
  Dans la suite on d\'efinira pour tout entier non nul $n$ pair 
   $\tilde J_{n}$ par
  $ \tilde J_{n}= J_{n}$.
  Enfin si $\alpha$ est un r\'eel positif de partie enti\`ere 
  $n_{\alpha}=2 p_{\alpha}$ on pose 
  $\tilde J_{\alpha}= \tilde J_{n_{\alpha}} \circ J_{\alpha-n_{\alpha}}$.
  \subsection{R\'esolution d'une \'equation diff\'erentielle factionnaire sur $]0,1[$ 
   sur $\alpha> 0$}
   Dans la suite on posera $D_{\alpha} (f) = D_{\alpha-n_{\alpha}}D_{n_{\alpha}}$
   pour tout r\'eel de partie enti\`ere paire $n_{\alpha}$.
Nous pouvons \'enonc\'e le th\'eor\`eme suivant :
   \begin{theoreme} \label{diffn}
  On consid\`ere un r\'eel positif $\alpha$, de partie enti\`ere $n_{\alpha}$.
 On suppose que $n_{\alpha} =2p_{\alpha}$ est un entier pair  et
  on se donne une fonction contractante $\psi$ localement contractante sur $]0,1[$ et dans 
  $L^{1}_{1-\alpha}([a,b])$. Alors 
  la fonction $\tilde J_{\alpha} (\psi)$ est solution sur $]0,1[$
  de l'\'equation diff\'erentielle 
  $D_{\alpha}(y)=\psi$
  avec les conditions initiales 
  $ y(0)=y^{(1)}(0)= \cdots=y^{(p_{\alpha}-1)}(0)=0$\\
  et 
  $ y(1)=y^{(1)}(1)= \cdots=y^{(p_{\alpha}-1)}(1)=0$.
 \end{theoreme}
 \begin{preuve}{}
Dans cette d\'emonstration nous poserons  $ \alpha'=\alpha-n_{\alpha}$.
Le cas $\alpha'=0$ ayant en fait \'et\'e trait\'e plus haut on supposera $0<\alpha'<1$.
La fonction $\psi$ \'etant contractante on a vu (voir le corollaire \ref{equdiff1}) que 
la fonction $y_{0} = J_{\alpha'}(\psi)$ est conctractante sur $[0,1]$ et 
v\'erifie, pour tout $x$ dans $]0,1[$, $ \left(D_{\alpha'}(y_{0})\right) (x) = \psi.$
On v\'erifie facilement que $y_{0}$ est continue sur $[0,1]$ la fonction $y=J_{n_{\alpha}}(y_{0})$ 
est d\'efinie sur $[0,1]$, et v\'erifie $D_{n_{\alpha}} (y) =y_{0}$ avec 
les conditions initiales  $ y(0)=y^{(1)}(0)= \cdots=y^{(p_{\alpha}-1)}(0)=0$\\
  et 
  $ y(1)=y^{(1)}(1)= \cdots=y^{(p_{\alpha}-1)}(1)=0$.\\
  On peut alors remarquer que par construction 
  $y= J_{n_{\alpha}}\circ J_{\alpha'} (\psi)=\tilde J_{\alpha} (\psi)$.
  On v\'erifie de plus, toujours par construction, pour tout $x$ dans $]0,1[$ 
 $ \left(D_{\alpha}(y)\right) (x)  = \left ( D_{\alpha'}\circ D_{n_{\alpha}}(y)\right)(x) 
 =  \left(D_{\alpha}(y_{0})\right) (x)= \psi (x)$,ce qui ach\`eve la d\'emonstration.
 \end{preuve}
  
       \subsection{Ecriture int\'egral de l'op\'erateur $\tilde {J_{\alpha}} f$}
      En tenant compte de la propri\'et\'e \ref{INTFRAC} nous pouvons \'ecrire 
      \begin{prop}
       Soit $\alpha = 2p_{\alpha}+\alpha'$ avec $p_{\alpha}\in \mathbb N^{\star}$
      et $0<\alpha'<1$.
      Pour une fonction $f$ v\'erifiant les hypoth\`eses du th\'eor\`eme
      \ref{contractante} et pour $\alpha>1$ non entier on a, pour tout r\'eel $x$ 
      dans $]0,1[$ : 
           $$\left(\tilde J_{\alpha}(f)\right)(x) =
      \frac{2^{\alpha}}{\Gamma (\alpha)} \int_{0}^{1} G_{p_{\alpha}} (x,y)
      \int_{0}^{y} (t-y)^{\alpha^{\prime}-1}f(t) dt dy.$$
            \end{prop}
      \section{D\'eriv\'ee sur un intervalle quelconque}
      \subsection{D\'eriv\'ee d'ordre $\alpha$, $0<\alpha<1$ sur un intervalle}
      Nous consid\'erons maintenant une fonction $f$ d\'efinie sur un intervalle 
      $[a,b], a<b$ quelconque. Nous noterons alors $f_{a,b}$ la fonction d\'efinie 
      pour tout $t\in [0,1]$ par $f(t) =  a + t (b-a)$. Si $\alpha\in [0,1]$ nous dirons que 
      $f$ est $\alpha$ fois d\'erivable en $x \in ]a,b[$ si $f_{a,b}$ est $\alpha$ d\'erivable en
      $\frac{x-a}{b-a}$.
      Nous poserons alors 
      $$\left( D_{\alpha,a,b} (f) \right) (x) = \left( D_{\alpha} (f_{a,b}) \right) \left( \frac{x-a}{b-a}\right).$$
      et, toujours pour une fonction $f$ d\'efinie sur $[a,b]$ nous d\'efinirons 
      $$\left( J_{\alpha,a,b} (f) \right) (x) = \left( J_{\alpha} (f_{a,b}) \right) \left(\frac{x-a}{b-a}\right),$$
      si, bien sr, cette derni\`ere quantit\'e est d\'efinie.
      Nous avons alors l'analogue de la propri\'et\'e \ref{invgauche}:
      \begin{prop} 
      Pour tout intervalle $[a,b]$, si $f$ et $\psi$ sont deux fonctions localement contractantes sur $]a,b[$ appartenant \`a 
      $L^1_{1-\alpha}([a,b]) $ telles que 
      $\psi =D_{\alpha,a,b} (f)$ on a pour tout 
      $x\in ]a,b[$
           $$ \left(J_{\alpha,a,b} (\psi)\right)(x)= f(x).$$
      \end{prop}
\begin{preuve}{}
Nous avons donc 
$ \left( D_{\alpha}(f_{a,b})\right) \left(\frac{x-a}{b-a}\right) 
= \psi (x)$ ce qui donne aussi, pour tout $t\in ]0,1[$ 
$ \left( D_{\alpha}(f_{a,b}) \right) (t) = \psi_{a,b} (t).$
C'est \`a dire que pour tout $t\in ]0,1[$ on a 
$ \left(J_{\alpha} (\psi_{a,b}) \right) (t) =  f_{a,b}(t) $
ou encore, en posant $t=\frac{x-a}{b-a}$
$  \left(J_{\alpha} (\psi_{a,b}) \right) (\frac{x-a}{b-a}) = 
 f_{a,b}(\frac{x-a}{b-a}) $,
 c'est \`a dire $ \left(J_{\alpha,a,b} (\psi)\right)(x)= f(x).$
 \end{preuve}
\subsection{Expression int\'egrale}
Dans ce paragraphe on suppose toujours $\alpha\in ]0,1[$.
On v\'erifie alors facilement les deux r\'esultats suivants :
\begin{prop}
Si $f$ est une fonction localement contractante sur  $]a,b[$ appartenant \`a $L^1_{1-\alpha}([a,b])$ nous pouvons \'ecrire pour tout $x\in ]a,b[$ 
$$
\left(D_{\alpha,a,b} (f)\right)(x)  =(b-a)^{\alpha}
  \frac{2^\alpha }{\Gamma (-\alpha)} 
\left(\int_{a}^x (x-u)^{-\alpha-1} \left ( f(u) - f (x)\right) du -
f(x) 
\frac{(x-a)^{-\alpha}}{\alpha}\right)
$$  
et 
$$
\left( J_{\alpha,a,b}(f)\right) (x) = \frac{(b-a)^{-\alpha}}{2^\alpha \Gamma (\alpha)}
\int_{a}^x \frac{f( u)}{(x-u)^{1-\alpha}}du.$$
\end{prop}

\subsection{D\'eriv\'ee enti\`ere sur un intervalle}
    Comme plus haut on pose pour tout entier naturel $n$ 
    $$  
   ( D_{n,a,b} f )(x) = \left(D_{n} (f_{a,b})\right) (\frac{x-a}{b-a}),
    $$
    et si $n$ pair on d\'efinit 
    $$ 
    (J_{n,a,b} f) (x)=\left(J_{n} (f_{a,b})\right) (\frac{x-a}{b-a})
    \quad \mathrm{si} \quad n=2p,
    $$
    $$ 
    (J_{n,a,b} f) (x)=\left(J_{1,a,b}\circ J_{n-1} (f_{a,b})\right) (\frac{x-a}{b-a})
    \quad \mathrm{si} \quad n=2p+1.
    $$
     On a comme plus haut, si $f$ est continue sur $[a,b]$
    $$ \left(D_{n,a,b} \circ J_{n,a,b}\right)(f) = f.$$
    On obtient aussi comme dans le cas des d\'eriv\'ees sur 
    $[a,b]$ d'ordre $\alpha$ avec $0<\alpha<1$
    \begin{prop}
    On suppose la fonction $f$ de classe $C^n$ sur $[a,b]$
    v\'erifiant les conditions initiales 
    $$ f^{(0)} (a) = f^{(1)} (a) =\cdots = f^{(p-1)} (a)=0,
     f^{(0)} (b) = f^{(1)} (b) =\cdots = f^{(p-1)} (b)=0.$$
    Alors, pour tout $x \in ]a,b[$ 
    $$\left(J_{2p,a,b} \circ D_{2p,a,b}(f) \right)  (x) = f (x).$$  
    \end{prop}
      \subsection{L'op\'erateur $\tilde J_{\alpha,a,b}$}
    On peut d\'efinir comme pr\'ec\'edemment l'op\'erateur $\tilde J_{\alpha,a,b}$ 
    \`a partir de $J_{n,a,b}$ pour $n$ entier pair et $J_{\alpha,a,b}$ pour $0<\alpha<1$. 
    C'est \`a dire que pour $\alpha>0$ de partie enti\`ere paire on pose 
    $\tilde J_{\alpha,a,b} = J_{n_{\alpha},a,b} \circ J_{\alpha-n_{\alpha},a,b}$
    avec $n_{\alpha} = [\alpha]$. Nous avons alors l'\'ecriture int\'egrale, si 
    $n_{\alpha}=2p_{\alpha}$, et pour $x\in ]a,b[$
  $$ \left(\tilde J_{\alpha,a,b}(f) (x)\right) = (-1)^{p }C^{-1}_{p_{\alpha},\alpha} 
    \int_{x}^{b} (x-a)^{p_{\alpha}}(u-a)^{p_{\alpha}} \Psi_{1,\alpha,a,b}(u,x)
    \Psi_{2,\alpha,a,b}(u,x) du$$
      avec $C_{p_{\alpha},\alpha} = (b-a)^{2p+\alpha} (p-1)!2^{\alpha}\Gamma(\alpha)$\\
    $$\Psi_{1,\alpha,a,b}(u,x) = \int_{u}^{b}\frac{(z-x)^{p_{\alpha}-1}  (z-u)^{p_{\alpha}-1} }
    {(z-a)^{2p_{\alpha}} }dz
   \quad \mathrm{et} \quad 
    \Psi_{2,\alpha,a,b}(u,x)= \int_{a}^{u}\frac{f(t)}{(u-t)^{1-\alpha}}dt .$$
        On obtient alors l' \'enonc\'e analogue au th\'eor\`eme \ref{GROSTHEO}
    \begin{theoreme}\label{GROSTHEOLOC}
  On consid\`ere un r\'eel positif $\alpha$, de partie enti\`ere $n_{\alpha}$.
 on suppose que $n_{\alpha} =2p_{\alpha}$ est un entier pair  et
  on se donne une fonction localement contractante $\psi$ sur $]a,b[$ et appartenant \`a $L^1_{1-\alpha}([a,b])$. Alors 
  la fonction $\tilde J_{\alpha,a,b} (\psi)$ est solution sur $]a,b[$
  de l'\'equation diff\'erentielle 
  $D_{\alpha,a,b}(y)=\psi$
  avec les conditions initiales 
  $ y(a)=y^{(1)}(a)= \cdots=y^{(p_{\alpha}-1)}(a)=0$\\
  et 
  $ y(b)=y^{(1)}(b)= \cdots=y^{(p_{\alpha}-1)}(b)=0$.
 \end{theoreme}
  \subsection{L'op\'erateur $D_{\alpha,\infty}$}
  \begin{remarque}
Si $f$ est une fonction contractante sur $[a,b]$, de support contenu dans $]a,b[$ 
nous pouvons\'ecrire, pour tout intervalle $[a',b']$ contenant $[a,b]$ 
et pour tout $x$ dans $]a,b[$ :
$$\frac{ \left(D_{\alpha,a,b} (f)\right)(x) }{ (b-a)^{\alpha}} -f(x) \frac{(x-a)^{-\alpha}}{\alpha}
= \frac{ \left(D_{\alpha,a',b'} (f)\right)(x) }{ (b'-a')^{\alpha}}-f(x) \frac{(x-a')^{-\alpha}}{\alpha}
$$ 
$$ J_{\alpha,a,b} (x) (b-a)^{\alpha} = J_{\alpha,a',b'} (x) (b'-a')^{\alpha}.$$ 
\end{remarque}
Cette conduite nous conduit \`a poser la d\'efinition 
\begin{definition}
On se donne un r\'eel $\alpha\in ]0,1[$.
On dit qu'une fonction $f$ est $D_{\alpha,\infty}$ d\'erivable en $x$ si 
$\displaystyle{ \lim_{a\rightarrow -\infty, b \rightarrow + \infty,a<x<b} 
\frac{\left(D_{\alpha,a,b} (f)\right)(x) }{(b-a)^{\alpha}}}$ existe et est fini.
On a alors 
$$\left( D_{\alpha,\infty} (f)\right) (x)=  \frac{2^\alpha }{\Gamma (-\alpha)} 
\int_{-\infty}^{x} (x-u)^{-\alpha-1} \left ( f(u) - f (x)\right) du.$$
\end{definition}
On pose alors 
$ \left( J_{\alpha,\infty} (f)\right) (x) = \displaystyle{ \lim_{a\rightarrow -\infty, b \rightarrow + \infty,a<x<b} 
\left(J_{\alpha,a,b} (f)\right)(x) (b-a)^{\alpha}}$ 
 si cette limite existe et est finie.
 Alors 
 $$ \left( J_{\alpha,\infty} (f)\right) (x) = 
 \frac{1}{2^\alpha \Gamma (\alpha)}
\int_{-\infty}^{x} \frac{f( u)}{(x-u)^{1-\alpha}}du.$$
\begin{remarque}
On remarque que si $f$ est une fonction \`a support compact les quantit\'es 
$\left( D_{\alpha,\infty} (f)\right) (x)$ et $ \left( J_{\alpha,\infty} (f)\right) (x)$ sont 
bien d\'efinies. 
\end{remarque}
\begin{remarque}
On peut remarquer que $J_{\alpha,\infty} $ correspond \`a la d\'eriv\'ee de Weyl-Liouville pour $0<\alpha<1$, et que $D_{\alpha,\infty}$ est une d\'eriv\'ee fractionnaire de Marchaud.
\end{remarque}
Nous pouvons alors \'enoncer la propri\'et\'e suivante :
\begin{prop}
Si $B$ est un r\'eel strictement positif fix\'e, et $\psi$ une fonction contractante sur $\mathbb R$ \`a support compact. Alors pour tout $x$ r\'eel la fonction $f(x) = (J_{\alpha,\infty}(\psi))(x)$ est solution de 
l'\'equation $D_{\alpha,\infty} (f) (x) = \psi(x).$
\end{prop}
\begin{preuve}{}
Les hypoth\`eses de la proposition impliquent que la quantit\'e 
$f(x) = (J_{\alpha,\infty}(\psi))(x)$ est d\'efinie pour tout r\'eel $x$.
Soit $x$ un r\'eel fix\'e.
Si $A$ est un r\'eel sup\'erieur \`a $x$ posons $f_{A} = 2A^{\alpha} J_{\alpha,-A,A} (\psi)$.
Alors d'apr\`es le th\'eor\`eme \ref{GROSTHEOLOC} tout 
$z\in]-A,A[$ v\'erifie $\frac{1}{(2A)^\alpha} \left( D_{\alpha,-A,A} (f_{A})\right) (z) = \psi (z).$
Soit $\epsilon>0$. Si $A$ est assez grand on peut trouver un r\'eel $B$, $A>B>0$,
$x\in ]-B,B]$ et  tel que pour $y\in [-B,B]$
$$ \vert \frac{1}{2A^{\alpha}} \left(D_{\alpha,-A,A} f_{A}\right) (y) 
-\left( D_{\alpha,\infty}f_{A}\right) (y) \vert \le \epsilon,$$
soit 
$$ \vert \psi (y) -\left(D_{\alpha,\infty} f_{A}\right) (y)\vert \le \epsilon.$$
D'autre part nous pouvons remarquer que pour tout $y\in [-B,B]$ 
$\ f(y) = f_{A}(y) $
 en posant $f = J_{\alpha,\infty} (\psi)$, 
et toujours en supposant  $A$ et $B$ suffisamment grands on conclut.
\end{preuve}
  \section{Appendice : calcul de $T_{2,k,l,R}.$}
  Pour un entier $k>$ on note $\gamma_{1,-k}$ et $\gamma_{2,-k}$ le coefficient de Fourier d'ordre $-k$ de $\left( \frac{1+R\chi}{1-R\bar \chi}\right)^{\alpha}$ 
  et de $\left( \frac{1- R \bar \chi}{1+R \chi}\right)^{\alpha}$.
  \subsection{ Calcul de $\gamma_{1,-k,R}$ et $\gamma_{2,k,R}$ pour $k$ suffisamment grand}
  Dans ce paragraphe on supposera 
  $ k =[Nx]$, avec $x>1$. Pour 
  $\alpha\in ]0,1[ $ on notera \'egalement 
  $\delta _{u,\alpha}$ le coefficient de Fourier de $\vert 1-\chi\vert ^{2\alpha}$. On sait que 
  $\delta_{u,\alpha}= C_{\alpha} \vert u\vert ^{-2\alpha-1}
  +o(u^{-2\alpha-1})$ avec $C_{\alpha}=  $.\\
  On peut remarquer que 
  $$ \left( \frac{1-R\chi}{1+R\bar \chi}\right)^{\alpha} 
  = \left(\frac{ \vert 1- R \bar \chi\vert ^2}{1-R^2\bar \chi^2}
  \right)^\alpha.$$
  On a donc :
  $$\gamma_{1,-k}=\Bigl  \langle \left(\sum_{v\in \mathbb Z}
  \delta_{v,\alpha} (R \bar \chi)^v \right) \left(\sum_{u\ge 0} 
  \beta_{u}^{(\alpha)} (R \bar \chi)^{2u}\right) \Big\vert
  \overline {\chi^k}\Bigr \rangle,$$
  soit 
$$
  \gamma_{1,-k}= \sum_{u \ge 0} \beta_{u} ^{(\alpha)}\delta_{k-2u,\alpha}
   R^k = \sum_{i=1}^4 S_{i}$$
  avec 
 $$ S_{1} =  \sum_{u = 0}^{k_{0}} \beta_{u} ^{(\alpha)} \delta_{k-2u,\alpha}
  R^k, \quad S_{2} = \sum_{u =k_{0}+1}^{(k-k_{0})/2-1}
  (-1)^{k-2u} \beta_{u} ^{(\alpha)}\delta_{k-2u,\alpha}
  R^k$$
  $$ S_{3}= \sum_{u= k-k_{0}/2} ^{(k+k_{0})/2} \beta_{u} ^{(\alpha)} \delta_{k-2u,\alpha} R^k, \quad S_{4}= \sum_{k+k_{0})/2+1} 
\beta_{u} ^{(\alpha)}\delta_{k-2u,\alpha} R^k,$$
  o\`u $k_{0}= N^\beta$ avec $\beta$ un r\'eel dans $]0,1[$ 
  qui sera fix\'e ult\'erieurement. Nous allons \'evaluer successivement les sommes $S_{1}, S_{2}, S_{3}, S_{4}$.
 Nous avons les estimations suivantes :
 \begin{enumerate} 
 \item
  \begin{equation} \label{eq1}
  \vert S_{1} \vert = O (k^{-2\alpha-1} (N^\beta)^\alpha 
  = O( N^{-2\alpha -1+ \beta \alpha}),
  \end{equation}
  \item
 $$
  \vert S_{2}\vert \le \frac{C_{\alpha}}{\Gamma (\alpha)} 
  \sum_{u =k_{0}+1}^{(k-k_{0})/2-1} u^{\alpha-1} (k-2u)^{-2\alpha-1} R^k \left (1+o(1)\right)$$
  d'o\`u l'existence d'un r\'eel $K>0$ tel que 
$$
   \vert  S_{2}\vert \le K \frac{C_{\alpha}}{\Gamma (\alpha)} 
    k_{0}^{\alpha-1} \sum_{u =k_{0}+1}^{(k-k_{0})/2-1}  (k-2u)^{-2\alpha-1} R^k 
    $$
    et donc 
    \begin{equation}\label{eq2}
   \vert  S_{2}\vert  =O( k_{0}^{-\alpha-1})
    \end{equation}
    \item
    \begin{align*}
    S_{3} &= \frac{1}{\Gamma (\alpha)} 
    \sum_{u= (k-k_{0})/2}^{(k+k_{0})/2} u^{\alpha-1}
    \delta_{k-2u} R^{k-2u} \left( 1+o(1)\right)\\
    &= -\frac{k^{\alpha-1} }{\Gamma (\alpha)} 
    \left( \sum_{u\le (k-k_{0})/2}
    \vert k-2u\vert ^{-2\alpha-1} R^{k-2u}
    \sum_{u\ge (k+k_{0})/2}
    \vert k-2u\vert ^{-2\alpha-1}R^{k-2u} \right) \left( 1+o(1)\right)          \end{align*}
     soit 
     \begin{equation}\label{eq3}
     \vert S_{3}\vert = O(k^{\alpha-1} k_{0}^{-2\alpha-1})=O(N^{\alpha-1-(2\alpha+1)\beta}).
     \end{equation}
     \item
   $$   \vert S_{4} \vert = O\left(
   \sum_{u= \frac{k+k_{0}}{2}} ^{+ \infty} 
   u^{\alpha-1} (2u- k)^{-2\alpha-1} \right)$$
   Il existe donc un r\'eel $K'$ strictement positif tel que  
   $$ 
   \vert S_{4}\vert \le K'
     \sum_{u= \frac{k+k_{0}}{2}} ^{+ \infty} 
   u^{\alpha-1} (2u- k)^{-2\alpha-1} R^k 
   $$
   Soit \begin{equation}\label{eq4}
    \vert S_{4}\vert = O ( N^{\alpha-1-2\alpha\beta}).
    \end{equation}
  \end{enumerate} 
  Les \'equations (\ref{eq1}), (\ref{eq2}), (\ref{eq3}), (\ref{eq4})
  permettent d'\'ecrire que pour un bon choix de 
  $\beta$ on a 
  \begin{equation}\label{eq5} \vert \gamma_{1,-k}\vert = 
  O (N^{-1-\alpha/2-\epsilon_{1}})
  \end{equation}
  avec $\epsilon_{1}>0$.\\
  
  De m\^{e}me on obtient 
  $$ \left( \frac{1-R\bar \chi}{1+R \chi}\right)^{\alpha} 
  = \left(\frac{ \vert 1- R  \chi\vert ^2}{1-R^2 \chi^2}
  \right)^\alpha.$$
  On a donc :
  $$\gamma_{2,k}=\Bigl  \langle \left(\sum_{v\in \mathbb Z}
  \delta_{v,\alpha} (R  \chi)^v \right) \left(\sum_{u\ge 0} 
  \beta_{u}^{(\alpha)} (R \chi)^{2u}\right) \Big\vert
 \chi^k\Bigr \rangle,$$
  c'est \`a dire qu'on obtient comme pr\'ec\'edemment 
  \begin{equation}\label{eq6} \vert \gamma_{2,k}\vert = 
  O (N^{-1-\alpha/2-\epsilon_{2}})
  \end{equation}
  avec $\epsilon_{2}>0$\\

  \subsection{Calcul des coefficients 
 $ \pi_{+}\left( \pi_{+}(\frac{\chi^{l} }{g_{2,R}}) \tilde \Phi_{N,R}\right)$ et 
  $ \pi_{+}\left( \pi_{+}(\frac{\chi^{k} }{\bar g_{1,R}}) \bar \Phi_{N,R}\right)$}
  Avec les notations introduites ci-dessus nous avons : 
  \begin {align*}
  \pi_{+}\left( \pi_{+}(\frac{\chi^{l} }{g_{2,R}}) \tilde \Phi_{N,R}\right) 
  &= \sum_{u=0}^{l}\beta_{l-u}^{(\alpha)} R^{l-u} \pi_{+}\left( \sum_{v\in \mathbb Z}
  \chi^{v+u-N-1}\gamma_{2,v}\right)\\
  &= \sum_{u=0}^{l}\beta_{l-u}^{(\alpha)} R^{l-u} \sum_{w=0}^{+\infty}  
  \chi^{w}\gamma_{2,w+N+1-u}
  \end{align*}
  et de m\^{e}me nous obtenons
  $$  \pi_{+}\left( \pi_{+}(\frac{\chi^{k} }{\bar g_{1,R}}) \bar \Phi_{N,R}\right) 
  = \sum_{v=0}^{k} \overline{ \beta_{k-v}^{(\alpha)}} R^{k-v}(-1)^{k-v} \sum_{w=0}^{+\infty} 
  \bar \gamma_{1,-(w+N+1-v)}\chi^{w}.$$
  \subsection{ Calcul de $ \left( I-H_{\tilde \Phi_{N}} H_{\Phi_{N}} \right)^{-1}  \pi_{+}\left( \pi_{+}(\frac{\chi^{l} }{g_{2,R}}) \tilde \Phi_{N,R}\right)$.}
  On a 
 \begin{align*}
 H_{\Phi_{N}}\left( \pi_{+}\left( \pi_{+}(\frac{\chi^{l} }{g_{2,R}}) \tilde \Phi_{N,R}\right)\right))  &=
 \sum_{u=0}^{l} \beta^{(\alpha)}_{l-u} R^{l-u} 
 \sum_{w_{0}=0}^{+\infty} \gamma_{2,w_{0}+N+1-u,R} \pi_{-}
 \left (  \sum_{v\in \mathbb Z} \gamma_{1,v,R}\chi^{w_{0}+v+N+1} \right)\\
 &= \sum_{u=0}^{l} \beta^{(\alpha)}_{l-u} R^{l-u} 
\sum_{w_{1}=0}^{+\infty}
 \left(\sum_{w_{0}=0}^{+\infty}  \gamma_{1,-(w_{1}+w_{0}+N+1)} \gamma_{2,w_{0}+N+1-u,R} \right)
 \chi^{-w_{1}}
 \end{align*}
  Si maintenant nous posons $ F_{1}(u,w_{1}) = 
 \displaystyle{  \sum_{w_{0}=0}^{+\infty}  \gamma_{1,-(w_{1}+w_{0}+N+1)} \gamma_{2,w_{0}+N+1-u,R}}$
  nous avons :
  \begin{align*}
  H_{\tilde \Phi_{N}} \left( H_{\Phi_{N}}\left(\pi_{+}\left( \pi_{+}(\frac{\chi^{l} }{g_{2,R}}) \tilde \Phi_{N,R}\right)\right)\right)&=  
   \sum_{u=0}^{l} \beta^{(\alpha)}_{l-u} R^{l-u} 
  \sum_{w_{1}=0}^{+\infty}F_{1}(u,w_{1}) \pi_{+} 
  \left(\sum_{v\in \mathbb Z}  \chi^{v-w_{1}-N-1} \gamma_{2,R,v}\right) \\
  &= \sum_{u=0}^{l} \beta^{(\alpha)}_{l-u} R^{l-u} 
 \sum_{w_{2}=0}^{+\infty} \left(\sum_{w_{1}=0}^{+\infty}F_{1}(u,w_{1})   
 \gamma_{2,w_{2}+w_{1}+N+1,R} \right) \chi^{w_{2}}
  \end{align*}
  ou encore, en posant $ F_{2}(u,w_{2}) = \displaystyle{ \sum_{w_{1}=0}^{+\infty}F_{1}(u,w_{1})   
 \gamma_{2,w_{2}+w_{1}+N+1,R}}$, 
 $$ H_{\tilde \Phi_{N}} \left( H_{\Phi_{N}}\left(\pi_{+}\left( \pi_{+}(\frac{\chi^{l} }{g_{2,R}}) \tilde \Phi_{N,R}\right)\right)\right)
 =  \sum_{u=0}^{l} \beta^{(\alpha)}_{l-u} R^{l-u} 
 \sum_{w_{2}=0}^{+\infty} F_{2}(u,w_{2}) \chi^{w_{2}}.
 $$ 
  On montre alors par r\'ecurrence que 
  $$  \left( I-H_{\tilde \Phi_{N}} H_{\Phi_{N}} \right)^{-1}  \pi_{+}\left( \pi_{+}(\frac{\chi^{l} }{g_{2,R}}) \tilde \Phi_{N,R}\right) = 
 \sum_{u=0}^{l} \beta^{(\alpha)}_{l-u} R^{l-u}  \sum_{m=0}^{+\infty} 
 \sum_{w_{2m}=0}^{+\infty} F_{2m}(u,w_{2m}) \chi^{w_{2m}}.
  $$
  Avec pour tout entier $n$, $2\le n\le m$
  $ F_{2n}(u,w_{2n}) = \displaystyle{\sum_{w_{2n-1}=0}^{\infty}F_{2n-1}(u,w_{p})   
 \gamma_{2,w_{2n}+w_{2n-1}+N+1,R}} $
 et 
 $ F_{2n-1}(u,w_{2n-1}) = 
  \displaystyle{ \sum_{w_{2n-2}}^{+\infty}  \gamma_{1,-(w_{2n-1}+w_{2n-2}+N+1)} }
  F_{2n-2}$ et avec 
  $F_{1},F_{2}$ comme plus haut.
  On obtient alors 
  $$T_{2,k,l,R} =  \sum_{u=0}^{l} \beta^{(\alpha)}_{l-u} R^{l-u}
   \sum_{v=0}^{k} \beta^{(\alpha)}_{k-v} (-1)^{k-v} R^{k-v} 
    \sum_{m=0}^{+\infty} 
 \sum_{w_{2m}=0}^{+\infty} F_{2m}(u,w_{2m})\bar \gamma_{1,-(w_{2m}+N+1-v),R}.
 $$
 Ce qui s'\'ecrit aussi, 
  \begin{align*}T_{2,k,l,R} &=  \sum_{u=0}^{l} \beta^{(\alpha)}_{l-u} R^{l-u}
   \sum_{v=0}^{k} \beta^{(\alpha)}_{k-v} R^{k-v} (-1)^{k-v}
    \sum_{m=0}^{+\infty} 
 \sum_{w_{2m}=0}^{+\infty} \bar \gamma_{1,-(w_{2m}+N+1-v),R}
 \\  &\sum_{w_{2n-1}=0}^{+\infty}  
 \gamma_{2,w_{2m}+w_{2m-1}+N+1,R} \cdots  \sum_{w_{0}=0}^{+\infty} 
  \gamma_{1,-(w_{1}+w_{0}+N+1)} \gamma_{2,w_{0}+N+1-u,R}.
  \end{align*} 
\subsection{Majoration de $T_{2,k,l,R}$}
Dans un premier temps nous supposons  $0<x<1, 0<y<1$. Posons $\epsilon = \min (\epsilon_{1}, \epsilon_{2})$ (avec les notations du paragraphe 4.3.1).
En utilisant les majorations 
$$ \vert   \gamma_{1,-(w_{1}+w_{0}+N+1)} \gamma_{2,w_{0}+N+1-u,R}\vert 
\le \frac{C_{1}}{(N+1-u)^{1+\alpha/2+\epsilon}}  
\frac{C_{1}}{(N+1+w_{0})^{1+\alpha/2+\epsilon}} $$
$$ \vert \gamma_{2,w_{2p}+w_{2p-1}N+1,R}\vert \le 
\frac{C_{1}}{(N+1+w_{2p-1})^{1+\alpha/2+\epsilon}} \quad \forall p\in [1,m]$$
$$ \vert \gamma_{1,- (w_{2p}+w_{2p+1}N+1),R}\vert \le 
\frac{C_{1}}{(N+1+w_{2p})^{1+\alpha/2+\epsilon}} \quad \forall p\in [1,m-1]$$
on obtient l'existence d' une constante $C$ telle que  pour $k,l$ v\'erifiant $k=[Nx]$, $l=[Ny]$ on ait 
$$ \vert T_{2,k,l,R}\vert \le
\left(  \sum_{u=0}^{l} \beta^{(\alpha)}_{l-u}\frac{C}{(N+1-u)^{1+\alpha/2+\epsilon}}\right)
\left( \sum_{v=0} ^{k}  \beta_{k-v}^{(\alpha)}  \frac{C}{(N+1-v)^{\alpha/2+\epsilon}}\right) 
\left( \sum_{m=0}^{+\infty} 
\left(\frac{C}{N^{\alpha/2+\epsilon}}\right)^{2m}\right),
$$
ce qui implique qu'il existe une constante $M$ telle que 
 \begin{equation} \label {eq7}
 \vert T_{2,k,l R}\vert \le M  \left( \sum_{u=0}^{l} \beta_{l-u}^{(\alpha)} 
 \frac{N^{-1-\alpha/2-\epsilon}}{(1-y)^{1+\alpha/2+\epsilon}}\right) \left( \sum_{v=0}^{k}\beta_{k-v}^{(\alpha)}  \frac{N^{-\alpha/2-\epsilon}}{(1-x)^{\alpha/2+\epsilon}}\right) 
 \end{equation}
 d'\`o\`u pour tout $0<\delta_{1}<\delta_{2}<1$ on a 
$ \vert T_{2,k,l R}\vert =o(N^{\alpha-1})$ uniform\'ement sur pour $x,y$ dans 
$[\delta _{1},\delta _{2}].$\\
Avec les notations ci-dessus supposons maintenant que 
 $ x\in [\delta _{1},\delta _{2}]$ et $y\in [\delta_{2},1]$.
 Nous pouvons 
 d\'ecomposer $T_{2,k,l,R}$ en deux sommes, 
 $\Sigma_{1}$ et $\Sigma_{2}$ avec 
   \begin{align*}\Sigma_{1} &=  \sum_{u=0}^{l-[N\delta_{1}]} \beta^{(\alpha)}_{l-u} R^{l-u}
   \sum_{v=0}^{k} \beta^{(\alpha)}_{k-v} R^{k-v} 
    \sum_{m=0}^{+\infty} 
 \sum_{w_{2m}=0}^{+\infty} \bar \gamma_{1,-(w_{2m}+N+1-v),R}
 \\  &\sum_{w_{2n-1}=0}^{+\infty}  
 \gamma_{2,w_{2m}+w_{2m-1}+N+1,R} \cdots  \sum_{w_{0}=0}^{+\infty} 
  \gamma_{1,-(w_{1}+w_{0}+N+1)} \gamma_{2,w_{0}+N+1-u,R}.
  \end{align*} 
  et 
    \begin{align*}\Sigma_{2} &=  
  \sum_{u=l-[N\delta_{2}]+1}^{l} \beta^{(\alpha)}_{l-u} R^{l-u}
   \sum_{v=0}^{k} \beta^{(\alpha)}_{k-v} R^{k-v} 
    \sum_{m=0}^{+\infty} 
 \sum_{w_{2m}=0}^{+\infty} \bar \gamma_{1,-(w_{2m}+N+1-v),R}
 \\  &\sum_{w_{2n-1}=0}^{+\infty}  
 \gamma_{2,w_{2m}+w_{2m-1}+N+1,R} \cdots  \sum_{w_{0}=0}^{+\infty} 
  \gamma_{1,-(w_{1}+w_{0}+N+1)} \gamma_{2,w_{0}+N+1-u,R}.
  \end{align*} 
On a comme pr\'ec\'edemment 
$\Sigma_{1} = O(N^{\alpha-1-\epsilon})$.
Par contre on peut \'ecrire, pour 
$ u\in \Bigl [  [N\delta_{2}]+1, N \Bigr]$
\begin{align*} \sum_{w_{0}=0}^{+\infty} 
  \gamma_{1,-(w_{1}+w_{0}+N+1)} \gamma_{2,w_{0}+N+1-u,R}
&\le \sum_{w_{0}=0}^{+\infty} 
\frac{C_{2}}{(N+1)^{1+\frac{\alpha}{2}+\epsilon} } \frac{C_{2}}{(w_{0}+1)^{1+\epsilon}}\\
&\le A_{\epsilon} \frac{C_{2}}{(N+1)^{1+\frac{\alpha}{2}+\epsilon} }
\end{align*}
avec 
$ A_{\epsilon} = \displaystyle {  \sum_{w_{0}=0}^{+\infty} \frac{C}{(w_{0}+1)^{1+\epsilon}}} $.
On obtient alors un r\'eel $C'$tel que 
$$
\Sigma_{2}
 \le
\left( \sum_{u=0}^{l} \beta^{(\alpha)}_{l-u}\frac{C'}{(N+1)^{1+\alpha/2+\epsilon}}\right)
\left( \sum_{v=0} ^{k}  \beta_{k-v}^{(\alpha)}  \frac{C'}{(N+1-v)^{\alpha/2+\epsilon}}\right) 
\left( \sum_{m=0}^{+\infty} \left(\frac{A C'}{N^{\alpha/2+\epsilon}}\right)^{2m}\right),
$$
et on retrouve 
 $\Sigma_{2}= o(N^{\alpha-1})$ et donc 
 $T_{2,k,l,R} =o(N^{\alpha-1})$ uniform\'ement pour 
tout $(x,y)$ dans $[\delta_{1}, \delta_{2}]\times [\delta_{2},1]$.\\
Supposons maintenant $y\in [0,\delta_{1}]$.
On a toujours l'\'equation (\ref{eq7}) et 
$ \displaystyle{\sum_{u=0}^l  \vert \beta_{l-u}^{(\alpha)}\vert } 
\le K_{\delta_{1}} N^\alpha$, o\`u $K_{\delta_{1}}$ est une constant qui ne d\'epend que de $\delta_{1}$.
On obtient donc encore une fois 
$T_{2,k,l,R} =o(N^{\alpha-1})$ uniform\'ement pour 
tout $(x,y)$ dans $[\delta_{1}, \delta_{2}]\times [0,\delta_{1}]$.\\
Nous pouvons donc \'enoncer la proposition \ref{propo3}

         \bibliography{Toeplitzdeux}

\end{document}